\documentclass[10pt]{article}

\usepackage{grffile}
\usepackage{graphicx}
\usepackage{epstopdf}

\usepackage{amsmath,amsfonts,amssymb}
\usepackage{amsthm}
\usepackage{color}
\usepackage[letterpaper,textheight=9in,textwidth=7in]{geometry}

\usepackage[skip=10pt,font=footnotesize]{caption}
\graphicspath{{}}

\newtheorem{Lemma}{Lemma}[section]

\newtheorem{Proposition}[Lemma]{Proposition}

\newtheorem{Remark}[Lemma]{Remark}
\newtheorem{Example}[Lemma]{Example}
\newtheorem{Definition}[Lemma]{Definition}

\newenvironment{Proof}[1][.]%
{\begin{trivlist}\item[]\textbf{Proof#1 }}%
{\qed\end{trivlist}}

 {\begin{trivlist}\item[]\textbf{Acknowledgments.}}{\end{trivlist}}

 {\begin{center}\textbf{Abstract}}{\end{center}}

\makeatletter\@addtoreset{figure}{section}\makeatother

\makeatletter \@addtoreset{equation}{section} \makeatother

\def\XXint#1#2#3{{\setbox0=\hbox{$#1{#2#3}{\int}$ }
\vcenter{\hbox{$#2#3$ }}\kern-.6\wd0}}

\newcommand{\R}{\mathbb{R}}
\newcommand{\C}{\mathbb{C}}
\newcommand{\N}{\mathbb{N}}

\newcommand{\ba}{\begin{align}}
\newcommand{\ea}{\end{align}}

\newcommand{\rmi}{\mathrm{i}}

\newcommand{\rmd}{\mathrm{d}}

\newcommand{\rme}{\mathrm{e}}
\newcommand{\rmo}{{\scriptstyle\mathcal{O}}}
\newcommand{\rmO}{\mathcal{O}}
\renewcommand{\Re}{\mathrm{Re}\,}

\newcommand{\eps}{{\varepsilon}}

\newcommand{\spec}{\mathrm{spec}\,}

\begin{document}
% \begin{center}
% {\fontsize{15}{15}\fontfamily{cmr}\fontseries{b}\selectfont{Signaling gradients in surface dynamics as basis for planarian regeneration}}\\[0.2in]
% Arnd Scheel$\,^1$, and Angela Stevens$\,^2$, and Christoph Tenbrock$\,^3$, \\[0.1in]
% \textit{\footnotesize 
% $\,^1$ University of Minnesota, School of Mathematics, 206 Church St. S.E., Minneapolis, MN 55455, USA\\
% $\,^2$University of M\"unster (WWU), Applied Mathematics, Einsteinstr. 62,
% D-48149 M\"unster, Germany \\
% $\,^3$ {{Munich, Germany, christoph.tenbrock@gmail.com }}
% }
% \date{\small \today} 
% \end{center}
\begin{center}
{\fontsize{15}{15}\fontfamily{cmr}\fontseries{b}\selectfont{Nonlinear eigenvalue methods for linear pointwise stability of nonlinear waves}}\\[0.2in]
Arnd Scheel \\[0.1in]
\textit{\footnotesize 
University of Minnesota, School of Mathematics,   206 Church St. S.E., Minneapolis, MN 55455, USA}
\date{\small \today} 
\end{center}

\begin{abstract}
We propose an iterative method to find pointwise exponential growth rates in linear problems posed on essentially one-dimensional domains. Such pointwise growth rates capture pointwise stability and instability in extended systems and arise as spectral values of a family of matrices that depends on a spectral parameter, obtained via a scattering-type problem. Different from methods in the literature that rely on computing determinants of this nonlinear matrix pencil, we propose and analyze an inverse power method that allows one to locate robustly the closest spectral value to a given reference point in the complex plane. The method finds branch points, eigenvalues, and resonance poles without a priori knowledge.
\end{abstract}

\setlength{\parskip}{4pt}
\setlength{\parindent}{0pt}

\section{Introduction}
Studying stability and instability of nonlinear waves and coherent structures informs our understanding of spatially extended nonlinear systems, with examples of applications that are of particular relevance to the present work ranging from instability in fluids \cite{doi:10.1146/annurev.fluid.37.061903.175810},  spatial ecology \cite{doi:10.1073/pnas.1420171112}, and biology \cite{PhysRevE.105.014602}, to material science \cite{PhysRevB.83.064113}. In models one analyzes stability of coherent structures using a variety of methods: explicitly \cite{MR1897705,MR1177566}, perturbatively \cite{MR1878337}, based on topological arguments \cite{MR3789546}, or, most often, using numerical methods that approximate the infinite domains by finite-domain  boundary-value problems \cite{stablab}. The analysis is commonly split into two parts, separating the stability in the far-field, with typically simple, spatially constant or periodic states, and the core region. The far field is usually more easily tractable, while detailed information on the core is rarely available explicitly or even asymptotically. In function spaces, the distinction between core and far-field is reflected in the distinction between point and essential spectra of the linearization, respectively; see \cite{fiedler03,MR3100266,sandstede02} for an overview and references therein. Essential spectra can be determined by algebraic computations after Fourier transform (or by solving boundary value problems after Bloch wave transforms in the case of asymptotically periodic states). Point spectra can be well approximated by problems in bounded domains with exponential convergence away from absolute spectra \cite{ssabs}.

Our focus here is on essentially one-dimensional systems, with one unbounded spatial direction, where spatial-dynamics methods have helped establish a wealth of results on existence and stability. Our interest is in identifying pointwise temporal growth rates, that is, exponential growth rates in time when initial conditions are compactly supported and growth is measured in a bounded region of space. One finds that such growth rates correspond to singularities in the spectral parameter $\lambda$ of the resolvent Green's function $G_\lambda(x,y)$  and we refer to those here as \emph{pointwise spectral values}. Such pointwise spectral values can not generally be identified as eigenvalues in an appropriate function space: they include resonances, that is, eigenvalues hidden by the essential spectrum, and branch points of the dispersion relation. Also, perturbation results for pointwise spectral values are more subtle: unlike spectra, they are in general not upper semicontinuous with respect to system parameters. 

Nevertheless, we propose here an iterative method that identifies pointwise spectral values using methods very much inspired by the power method, which is at the heart of computational methods for most eigenvalue problems. As a specific objective, we focus on a basic algorithmic challenge: given a reference point $\lambda_0\in\C$:
\begin{center}
 \emph{Find the pointwise spectral value $\lambda$ closest to $\lambda_0$!}
\end{center}
% 
% \begin{enumerate}
%  \item find the closest pointwise spectral value in constant-coefficient problems;
%  \item find the closest pointwise spectral value in variable coefficient problems, allowing for instance embedded eigenvalues and resonances.
%  \end{enumerate}
Questions of this type arise when investigating  resonances in Schr\"odinger operators and in nonlinear optics, although algorithms of the nature proposed here do not appear to have been used in the literature. Even in constant- or periodic-coefficient problems, such tasks present challenging problems, relating to many questions in fluid mechanics  \cite{doi:10.1146/annurev.fluid.37.061903.175810,vansaarloos03}, material science \cite{PhysRevB.83.064113}, and ecology \cite{doi:10.1073/pnas.1420171112}.  Current methods require an intricate parameter continuation of eigenvalue problems and may at times miss leading pointwise growth rates; see for instance \cite{brevdo_linear_1999,MR2183609}. 

Our focus on pointwise spectral values originates in work on pointwise Green's functions in the context of shock stability  \cite{zumbrun98}. We are further motivated by the inherently pointwise nature of the analysis of coherent structures and the Evans function in many examples \cite{MR1068805}, the vast literature in fluid dynamics concerned with convective and absolute (pointwise) instabilities \cite{doi:10.1146/annurev.fluid.37.061903.175810}, and, lastly, the role of pointwise stability in the selection of fronts propagating into unstable states \cite{holzerscheel14,https://doi.org/10.48550/arxiv.2012.06443}. % avery scheel, holzer faye siemer scheel
Our point of view is shaped by the perspective of \emph{nonlinear eigenvalue problems}, that is, matrix or operator families that depend nonlinearly on a spectral parameter and where spectral parameter values for which the inverse of the operator is not analytic
are the object of interest. This point of view allows us to simultaneously treat far-field and core, to preserve structure of eigenvalue problems, and to develop iterative methods that provably converge to leading eigenvalues. Theoretically, our first contribution is a formulation of the problem of finding pointwise spectral values as a nonlinear eigenvalue problem, where local power series are readily computed from a homological equation. Our second contribution develops an inverse power method for this nonlinear eigenvalue problem that provably converges to the nearest spectral value. We prove in particular that, curiously, the method detects eigenvalues even past the radius of convergence of the local power series expansion.

The approach developed  here is complementary to Evans function methods.  The Evans function is a  popular and well-developed analytical and computational tool for the analysis of point spectra and resonances, a Wronskian-type complex analytic function that enables one to find eigenvalues as roots of an analytic function, exploiting for instance winding number computations to count numbers of unstable eigenvalues and to thereby establish robustly stability or instability; see for instance \cite{MR1068805,sandstede02}. The Evans function is computed either via differential forms or, more directly, taking a determinant of bases of bounded solutions to the linearized equation at spatial $\pm \infty$. It can in fact be related to an operator-theoretic, non-pointwise Fredholm determinant \cite{MR2350362}. The in many ways most challenging problems arise when studying point spectra located near or embedded in essential spectra. The approach here provides a more canonical computational view on these spectral problems while, at the same time, emphasizing the {pointwise} character of the stability questions of interest. By avoiding determinants, it has potential to perform better in large systems. 

% In this paper, our emphasis is on the theoretical aspect of this approach and on demonstrating the effectiveness in a variety of test examples, hoping that future work will demonstrate the viability in more complex situations. 

% 
% 
% In both problems the term ``pointwise spectral value'' refers to the closest singularities of a \emph{pointwise resolvent} to a linear operator, that is, the closest values of the spectral parameter $\lambda$ where the Green's function $G\lambda(x,y;\lambda)$ is not analytic for fixed $x$ and $y$. Those values are the effective obstructions to deformations of contour integrals in the complex plane when constructing heat kernels via inverse Laplace transform and thereby determine pointwise exponential growth rates; see for instance \cite{holzerscheel14,zumbrun98}.

\paragraph{Outline.} The remainder of the paper is organized as follows. We set up a somewhat general framework for eigenvalue problems and formulate the nonlinear pointwise eigenvalue problem in \S\ref{s:2}. We discuss an inverse power method for nonlinear eigenvalue problems and its convergence properties in \S\ref{s:3}, and discuss implementation, both for the inverse power method and for the derivation of the nonlinear eigenvalue problem on the Grassmannian, in  \S\ref{s:4}. We conclude with example computations of pointwise spectral values in constant and variable-coefficient problems in \S\ref{s:5} and a brief summary in \S\ref{s:6}.

\section{Pointwise nonlinear eigenvalue problems from linearization at heteroclinic profiles}\label{s:2}

\subsection{First-order ODEs from eigenvalue problems}
We consider eigenvalues problems that arise in the linearization at traveling waves, of the form
\begin{equation}\label{e:twlin}
 u_x=A(x;\lambda)u,\quad x\in\R,\ u\in\C^N,
\end{equation}
with matrix coefficients $A(x;\lambda)\in \C^{N\times N}$, continuous in $x$ and analytic in $\lambda$. We focus on the simplest case of asymptotically constant coefficients
\begin{equation}\label{e:aprop}
 \lim_{x\to\pm\infty}A(x;\lambda)=A_\pm(\lambda).
\end{equation}
These equations arise when casting the linearization in the comoving frame as a first-order ODE, substituting $\rme^{\lambda t}$ for time dependence. 
\begin{Example}\label{ex:1} We explain the transformations in the case of a simple example, the scalar nonlinear diffusion equation 
\begin{equation}\label{e:kpp}
 w_t=w_{xx}+w-w^3,
\end{equation}
with traveling fronts $w=w_*(x-ct)$ connecting $w=w_-$ at $x=-\infty$ to $w=w_+$  at $x=+\infty$, $w_\pm\in \{-1,0,1\}$. The linearization at such a front satisfies
\begin{equation}\label{e:kppl}
 w_t=w_{xx}+cw_x + (1-3w_*^2)w =: \mathcal{L}w,
\end{equation}
which leads to the formulation in the form \eqref{e:twlin},
\begin{equation}
 u_x=A(x;\lambda)u,\qquad A(x;\lambda)=
 \begin{pmatrix}
 0&1\\
 -1+3w_*^2(x)+\lambda & -c
 \end{pmatrix},
 \end{equation} with 
 \begin{equation}
A_\pm(\lambda)=
 \begin{pmatrix}
 0&1\\
 -1+\lambda & -c
 \end{pmatrix}, \text{ if }w_\pm=0, \text{ or }
 \quad A_\pm(\lambda)=
 \begin{pmatrix}
 0&1\\
 2+\lambda & -c
 \end{pmatrix}, \text{ if }|w_\pm|=1.
\end{equation}
\end{Example}
Such a formalism has been extended to many other situations, including asymptotically periodic coefficients $A_\pm=A_\pm(x;\lambda)=A_\pm(x+L_\pm;\lambda)$ or ill-posed equations on an infinite-dimensional state space $u\in X$ for problems in infinite cylinders or modulated waves and it would be interesting to pursue the methods developed here in such contexts as well \cite{ssmodstab,ssmorse,MR1759902}. We note that we explicitly allow nonlinear, polynomial dependence of $A(x;\lambda)$ on $\lambda$, for cases with higher-order time derivatives, for instance the wave equation, or for cases where the spectral parameter is replaced by a polynomial to resolve branch points in the dispersion relation; see for instance Examples  \ref{ex:source} and \ref{ex:1ctd}, below. 

% COMMENTS: STRUCTURE OF EVP: MULTISYMPLECTIC STUFF, Maslov, Krein, symmetries such as reversibilities

One can in much generality relate properties of the operator $\mathcal{T}(\lambda)=\frac{\rmd}{\rmd x}-A(x;\lambda)$ to properties of the linearization of the traveling wave, in our example the operator $\mathcal{L}$, both in function spaces and in a pointwise sense; see for instance \cite{sandstede02,ssmorse,holzerscheel14}.
We will therefore focus on properties of the (linear) operator pencil $\mathcal{T}$ without trying to relate back to the traveling-wave linearization in any generality. 

It is not hard to see \cite{palmer,ssmorse} that $\mathcal{T}(\lambda)$ is Fredholm as a closed, densely defined operator on, say, $L^2(\R,\C^N)$ with domain of definition $H^1(\R,\C^N)$ if and only if the asymptotic matrices $A_\pm(\lambda)$ are hyperbolic, that is, $\spec A_\pm(\lambda)\cap \rmi\R = \emptyset$. The Fredholm index is then given by the difference of Morse indices,
\begin{equation}\label{e:fm}
 \mathrm{ind}\,(T(\lambda))=i_\mathrm{M}(A_-(\lambda))-i_\mathrm{M}(A_+(\lambda)),
\end{equation}
where $i_\mathrm{M}(A)$ counts the eigenvalues of $A$ with positive real part with multiplicity; see for instance \cite{ssmorse} and references therein. For well-posed equations, $\mathcal{L}-\lambda$ and thereby $\mathcal{T}(\lambda)$ are invertible for $\Re\lambda\gg 1$, such that the Morse index there is constant, $i_\mathrm{M}(A_+(\lambda)\equiv i_\infty=i_\mathrm{M}(A_-(\lambda)$.  Fredholm properties, that is, closedness of range and dimensions of kernel and cokernel,  of $\mathcal{T}(\lambda)$ and of $\mathcal{L}-\lambda$ agree. 

In the Fredholm 0 region, the analytic Fredholm theorem guarantees that generalized multiplicities of isolated eigenvalues of $\mathcal{L}$ are finite. In fact, generalized multiplicities of an eigenvalue $\lambda$ of  $\mathcal{L}$ agree with the multiplicity of an eigenvalue of $\mathcal{T}(\lambda)$ when the latter is defined  as follows; see \cite{trofimov,mennicken,gohberg} for the introduction of this concept and context, respectively.
\begin{Definition}[Algebraic multiplicities and Jordan chains]\label{d:jch}
Suppose $\mathcal{T}(\lambda_*)$ is Fredholm of index 0 with nontrivial kernel. We say a polynomial $u(\lambda)$ of order $p$, is a root function if $T(\lambda)u(\lambda)=\rmO((\lambda-\lambda_*)^{p+1})$. For root functions $u(\lambda)=\sum_{j=0}^p u_j(\lambda-\lambda_*)^j$, we refer to the $u_j$, $j<p$ as generalized eigenvectors. Note that $u_p$ is always an eigenvector, that is, $T(\lambda_*)u_p=0$. We define the algebraic multiplicity of $\lambda_*$ as the dimension of the (linear) space of root functions (of arbitrary degree $p$). 
\end{Definition}
A quick calculation verifies that the definitions here agree with the usual definitions of algebraic multiplicity in the case of standard eigenvalue problems. 

\begin{Example}
 In our example, a generalized eigenvector to $\lambda=0$ of $\mathcal{L}$ solves $\mathcal{L}w_1+w_0=0$, $\mathcal{L}w_0=0$. Defining $u_j=(w_j,w_{j,x})$, $j=0,1$, we find immediately from algebraic manipulation that $\mathcal{T}(0)u_0=0$ and  $\mathcal{T}(0)u_1+\mathcal{T}'(0)u_0=0$, showing how  Jordan chains are equivalent. 
\end{Example}

Since we did not formally introduce a general class of operators $\mathcal{L}$, we only state informally that in addition to Fredholm properties, also algebraic multiplicities of eigenvalues in the Fredholm index 0 region coincide for $\mathcal{L}-\lambda $ and $\mathcal{T}(\lambda)$.

\subsection{The Grassmannian and pointwise formulations of eigenvalue problems}
Our aim here is to develop a pointwise-in-$x$ formulation of the spectral problem for $\mathcal{T}(\lambda)$. Such formulations have been used extensively in the context of Schr\"odinger operators and developed also more generally in connection with stability of nonlinear waves in \cite{zumbrun98}. 
We start by considering the ODE \eqref{e:twlin} in the Fredholm index 0 regime where $i_\mathrm{M}(A_\pm(\lambda))=i_\infty$. The linear equation induces  a flow on $k$-dimensional (complex) subspaces $\mathrm{Gr}(k,N)$. We write $E^\mathrm{s/u}_\pm(\lambda)$ as the generalized eigenspaces of $A_\pm(\lambda)$ to eigenvalues $\nu$ with $\Re\nu<0$ and $\Re\nu>0$, respectively. 
These subspaces are invariant under $A_\pm(\lambda)$, respectively, and thereby invariant under the flow to $u'=A_\pm(\lambda)u$. One finds that $E^\mathrm{s}_+(\lambda)$ is unstable and  $E^\mathrm{u}_-(\lambda)$ is stable for the dynamics on $\mathrm{Gr}(N-i_\infty,N)$ and $\mathrm{Gr}(i_\infty,N)$ , respectively, that is, eigenvalues of the linearization at those equilibria all have positive or negative real part, respectively. 
One can then find unique subspaces $E^\mathrm{s}_+(x;\lambda)$ and  $E^\mathrm{u}_-(x;\lambda)$, continuous in $x$ and locally analytic in $\lambda$,  
which are invariant under the flow on the Grassmannian induced by \eqref{e:twlin} and converge to $E^\mathrm{s}_+(\lambda)$ and  $E^\mathrm{u}_-(\lambda)$, 
for $x\to +\infty$ and $x\to -\infty$, respectively. In particular, $\lambda$ is an eigenvalue if and only if $E^\mathrm{s}_+(0;\lambda)\cap E^\mathrm{u}_-(0;\lambda)\neq \{0\}$ is nontrivial. 

\begin{Lemma}[Analytic bases]\label{l:ab}
 For any fixed compact region $\Omega\subset\C$ where $E^\mathrm{s/u}_\pm(0;\lambda)$ are analytic, there exist analytic bases  $w_j^\mathrm{u}(\lambda)$, $1\leq j\leq i_\mathrm{M}$ and $w_j^\mathrm{s}(\lambda)$, $i_\mathrm{M}+1\leq j\leq N$ that span $E^\mathrm{s/u}_\pm(0;\lambda)$, respectively. 
\end{Lemma}
\begin{Proof}
The existence of such bases is an immediate consequence of \cite[Rem. 2]{shubin}, which guarantees the existence of an analytic complement and thereby analytic projections onto $E^\mathrm{s/u}_\pm(0;\lambda)$, respectively, and \cite[22,\S II.4.2]{kato}, which concludes the existence of analytic bases for subspaces given as the range of an analytic projection. A more constructive approach was described in \cite{MR2221065}, constructing analytic bases to $E^\mathrm{s/u}_+(\lambda)$, first, lifting them to nearby subspaces at $x=\pm L$, $L\gg 1$, and then transporting bases with the flow to the ODE \eqref{e:twlin}.
 
We describe a third approach here that relates to our specific choice of bases, below. Write $E(\lambda)$ for an analytic family of subspaces, either $E^\mathrm{s}_+(0;\lambda)$ or $E^\mathrm{u}_-(0;\lambda)$,   choose a complement $F_0$ for $E_0:=E(\lambda_0)$, and choose a basis $w_1,\ldots,w_m$ in $E(\lambda_0)$. Write $P_0$ for the projection along $F_0$ onto $E(\lambda_0)$. The subspace $E(\lambda)$ is then given as the graph of a map $H(\lambda):E_0\to F_0$, whenever $E(\lambda)\cap F_0=\{0\}$. We claim that the coefficients of $H(\lambda)$ have isolated poles of finite order, only, whenever $E(\lambda)\cap F_0\neq\{0\}$. For this, fix $\lambda_1$ where $H(\lambda)$ is singular, and choose $E_1,F_1$ complementary subspaces so that $E(\lambda)=\mathrm{graph}\,(H_1(\lambda))$, $H_1(\lambda):E_1\to F_1$ analytic for $\lambda\sim \lambda_1$. The map $H(\lambda)$ is then explicitly found from $H(\lambda)=(1-P_0)(\mathrm{id}+H_1(\lambda))\left(P_0(\mathrm{id}+H_1(\lambda) \right)^{-1}$, where the  inverse yields a meromorphic function with isolated poles. 

We therefore find basis vectors $W_j(\lambda)=w_j+H(\lambda)w_j$, $1\leq j\leq m$, for all $\lambda$ except for a finite  set of points where the $W_j$ have poles. For each of the $W_j$, we can however remove the pole singularity at a point $\lambda_\ell$ multiplying the singular basis vector $W_j$ by $(\lambda-\lambda_\ell)^p$, where $p$ is the maximal order of the pole in the components of $W_j$. We thereby obtain analytic vectors $\tilde{W}_j$ which form a basis for all $\lambda$.

%  
% project with P, fixed, on a fixed reference space so rthat at some lambda_0 the projection is one-to-one
% 
% choose a fixed basis in the range of P and construct the basis in the family E(lambda) from P^{-1}
% 
% this involves solving a system of linear equations which can be achieved over the field of meromorphic functions
% 
% we thus obtain a basis  of vectors which are each meromorphic functions
% 
% in particular, they are linearly independent over the field of meromorphic functions
% 
% we now 'regularize the basis by multiplying each vector with (lambda-lambda_l)^k_l, where lambda_l are the poles and k_l the order of the poles. Since the basis was linearly independent over the meromorphic functions, this operation preseves linear independence and yields an analytic basis

\end{Proof}

The same result applies in the case where bases have branch points which are resolved writing $\lambda=\varphi(\gamma)$. Subspaces that are analytic in $\gamma$ then have analytic bases.

\begin{Definition}[Pointwise eigenvalue problem]\label{d:iota}
 We define the trivialization of the bundles $E^\mathrm{s}_+(0;\lambda)$ and $E^\mathrm{u}_-(0;\lambda)$ through maps 
 \begin{align*}
  \iota^\mathrm{u}(\lambda):\C^{i_M}\to E^\mathrm{u}_-(0;\lambda),\qquad & u\mapsto \sum_{j=1}^{i_\mathrm{M}}
  u_j w_j^\mathrm{u}(\lambda),\\
  \iota^\mathrm{s}(\lambda):\C^{N-i_M}\to E^\mathrm{s}_+(0;\lambda),\qquad & u\mapsto \sum_{j=i_\mathrm{M}+1}^{N}
  u_j w_j^\mathrm{s}(\lambda) ,  
 \end{align*}
where the bases $w_j^\mathrm{s/u}(\lambda)$ were constructed in Lemma \ref{l:ab}. We then define the intersection map 
\[
 \iota_\mathrm{sec}(\lambda):E^\mathrm{u}_-(0;\lambda)\times E^\mathrm{s}_+(0;\lambda)\to \C^N, (w^\mathrm{u},w^\mathrm{s})\to w^\mathrm{u}-w^\mathrm{s},
\]
and its trivialization 
\begin{equation}
 \iota(\lambda)=\iota_\mathrm{sec}(\lambda)\circ \left(\iota^\mathrm{u}(\lambda),\iota^\mathrm{s}(\lambda)\right).
\end{equation}
We also define the associated Evans function
 \begin{equation}\label{e:Evans}
  \mathcal{E}(\lambda)=\mathrm{det}\,\iota(\lambda).
 \end{equation}
\end{Definition}

\begin{Proposition}\label{e:ptwisemult}
 The nonlinear eigenvalue problems $\mathcal{T}(\lambda)$ and $\iota(\lambda)$ are equivalent in the sense that geometric and algebraic multiplicities, in a region $\Omega$ where $\mathcal{T}(\lambda)$ is Fredholm index 0. In particular, the algebraic multiplicity of eigenvalues of  $\mathcal{T}(\lambda)$ equals the order of the root of the Evans function $ \mathcal{E}(\lambda)=\mathrm{det}\,\iota(\lambda)$. 
\end{Proposition}
\begin{Proof}
We claim that root functions for $\mathcal{T}$ and $\iota$ are in 1-1 correspondence. Indeed, given a root function $u^0(\lambda)$ for $\iota$, we can construct functions $u(x;\lambda)$ by solving the initial-value problem at $x=0$ and find bounded solutions up to the order of the root function. Conversely, restricting root functions for $\mathcal{T}$ to $x=0$ yields root functions for $\iota$. For finite-dimensional nonlinear eigenvalue problems as the one defined by $\iota$, the algebraic multiplicity is as defined in Definition \ref{d:jch} and agrees with the order of the root of the determinant \cite{trofimov}.
\end{Proof}
We are also interested in a version of Proposition \ref{e:ptwisemult} concerned with the analytic extension of $\iota(\lambda)$ past the essential spectrum. As an analytic function, $\iota$ has a uniquely defined analytic extension to some open set $\Omega\subset \C$. The motivation for considering this extension is rooted in the relation between this extension of $\iota$ and pointwise singularities of the Green's function.

\begin{Proposition}[Singularities of the pointwise Green's functions and $\iota$]\label{p:ptwiseext}
 Consider the Green's function of $\mathcal{T}(\lambda)$, solution to $\mathcal{T}(\lambda)G(x,y;\lambda)=\delta(x-y)\mathrm{id}$. Then $G(x,y;\lambda)$ with $x,y$ fixed, arbitrary, possesses an analytic extension in $\lambda$ into the region where  $\iota(\lambda)^{-1}$ possesses an analytic extension. On the other hand, $G(x,y;\lambda)$ is not analytic when 
 \begin{enumerate}
 \item $E^\mathrm{u}_-(0;\lambda)$ or $E^\mathrm{s}_+(0;\lambda)$ are not analytic, or when 
 \item $E^\mathrm{u}_-(0;\lambda)$ and $E^\mathrm{s}_+(0;\lambda)$ intersect nontrivially. 
\end{enumerate}
\end{Proposition}
Note that the poles of $\iota(\lambda)$ do not necessarily contribute to singularities of $\iota(\lambda)^{-1}$. The case(ii) corresponds to zeros of an extension of the Evans function, yielding resonances or embedded eigenvalues, both of which we refer to as extended point spectrum, following \cite{ssabs,rademacher07}. Analyticity of  $E^\mathrm{u}_-(0;\lambda)$ and $E^\mathrm{s}_+(0;\lambda)$ follows from analyticity of  $E^\mathrm{u}_-(\lambda)$ and $E^\mathrm{s}_+(\lambda)$ with sufficiently rapid convergence of the matrices $A(x;\lambda)$ by  results usually referred to as ``Gap Lemmas'' \cite{kapitula98,gardner98}. Absent such conditions, subspaces  $E^\mathrm{u}_-(0;\lambda)$ and $E^\mathrm{s}_+(0;\lambda)$ may exhibit essential singularities \cite{sandstede04}. Singularities of the asymptotic subspaces correspond to branch point singularities at infinity, since subspaces are obtained from algebraic equations; see \cite{holzerscheel14} for an extensive discussion of those singularities, referred to there as right-sided pointwise growth modes. 

\begin{Proof}[ of Prop. \ref{p:ptwiseext}]
Setting without loss of generality $y=0$, we need to solve $\mathcal{T}(\lambda)G(x,0;\lambda)=\delta(x)v$, $v\in\C^N$. Clearly, this requires a solution to the ODE defined by $\mathcal{T}$ with a jump at $x=0$ of size $v$. In the region where $\mathcal{T}$ is invertible, such a solution can be obtained uniquely by solving $\iota(\lambda)(w^\mathrm{u},-w^\mathrm{s})=v$, and extending the initial condition $u_-=\sum_{j=1}^{i_\mathrm{M}} w^\mathrm{u}_j u_j$ to $x<0$ and extending the initial condition $u_+=\sum_{j=i_\mathrm{M}+1}^{N} w^\mathrm{s}_j u_j$ to $x>0$. This construction clearly shows analyticity of $G$ given analyticity of $\iota^{-1}$, and, on the other hand, that conditions (i) and (ii) are necessary for analyticity of $G$. 
\end{Proof}
Information on the Green's kernel $G$ translates via Laplace transform directly into pointwise information on solutions to $\rme^{\mathcal{L}t}$ which we state here only informally. Given compactly supported initial conditions $w_0(x)$, $\sup_{|y|\leq K}\left(\rme^{\mathcal{L}t}u_0\right)(y)$ decays uniformly for any $K $ if $\iota(\lambda)$ is analytic in $\{\Re\lambda\}>0$. Conversely, the supremum grows exponentially if $\iota(\lambda)$ has a singularity in $\{\Re\lambda>0\}$ since direct Laplace transform of the heat kernel would otherwise imply analyticity of $G$; see for instance \cite[Cor. 2.3]{holzerscheel14}.
In a way similar to the case of point spectrum, one can associate Jordan chains to points $\lambda$ where $\iota$ is not invertible. 
% 
% A natural question then is to determine the singularities of $\iota(\lambda)^{-1}$. We list here the three sources of singularities:
% \begin{enumerate}
%  \item $E^\mathrm{u}_-(\lambda)$ or $E^\mathrm{s}_+(\lambda)$ are not analytic;
%  \item  $E^\mathrm{u}_-(0;\lambda)$ or $E^\mathrm{s}_+(0;\lambda)$ are not analytic;
%  \item $E^\mathrm{u}_-(0;\lambda)$ and $E^\mathrm{s}_+(0;\lambda)$ intersect nontrivially. 
% \end{enumerate}
% The cases (ii) and (iii) do not depend on $x$ since subspaces can be pushed from $x_1$ to $x_2$ using the analytic solution operator to \eqref{e:twlin}. The case (i) corresponds to branch point singularities at infinity, since subspaces are obtained from algebraic equations; see \cite{holzerscheel14} for an extensive discussion of those singularities, referred to there as right-sided pointwise growth modes. The case (ii) can be excluded with sufficiently rapid convergence of the matrices $A(x;\lambda)$, results usually referred to as ``Gap Lemmas'' \cite{kapitula98,gardner98}. Absent such conditions, subspaces may exhibit essential singularities \cite{sandstede04}. Case (iii), corresponds to zeros of an extension of the Evans function, yielding resonances or embedded eigenvalues, both of which we refer to as extended point spectrum here following \cite{ssabs,rademacher07}. 
% 
% what really happens when convergence is not strong enough?? essential singularities? Can we fix a singularity at infinity

In the following, we assume that a meromorphic realization of $\iota$ via meromorphic choices of bases, that is, of trivializations $\iota^\mathrm{u/s}$, has been fixed in the region where $E^\mathrm{u}_-(0;\lambda)$ and $E^\mathrm{s}_+(0;\lambda)$ are analytic. 
\begin{Definition}[Spectral values]\label{d:sv}
 We say $\lambda_0$ is a spectral value of $\iota$ if $\iota^{-1}(\lambda)$ is not analytic at $\lambda_0$.  Equivalently, conditions (i) or (ii) in Proposition \ref{p:ptwiseext} are violated.  
\end{Definition}

\begin{Remark}[Removing branch points]
 Singularities stemming from singularities of the asymptotic subspaces are branch points and can be removed using a polynomial reparametrization of the spectral parameter, $\lambda=\varphi(\gamma)$. Considering the new spectral problem with eigenvalue parameter $\gamma$, all of the above considerations apply again.
\end{Remark}

% 
% 
% 
% how does this relate to singularities of $\iota^{-1}$
% 
% do we need to say something??? probably yes...
% 
% * given lambda, check first if eigenspaces are analytic 
% * note that non-analyticity can only arise through branch points...
% 
% * then introduce locally analytic bases
% * then check $\iota$ inverse, whose local singularities then do not depend on choice of bases
% 

% 
% 
% absolute spectrum? branch points
\begin{Example}\label{ex:source}
 As a simple first example, we consider 
 \[
  w_t=w_{xx}-2\,\mathrm{sign}(x) w_x,
 \]
which leads to the spatial ODE
\begin{equation}\label{e:source}
 u_x=v,\qquad v_x=2\,\mathrm{sign}(x) v+\lambda u,
\end{equation}
with 
 \[
E_+^\mathrm{s}(\lambda)=\begin{pmatrix}1\\1-\sqrt{1+\lambda} \end{pmatrix},\qquad 
E_-^\mathrm{u}(\lambda)=\begin{pmatrix}1\\-1+\sqrt{1+\lambda}\end{pmatrix},
\]
and 
\[
 \mathcal{E}(\lambda)=2\left( 1-\sqrt{1+\lambda}\right).
\]
We find a zero at $\lambda=0$, case (iii) above, and a branch point at $\lambda=-1$, case (i). Note that the branch point corresponds to a spectral value of $\iota$, which can be removed by passing to a Riemann surface, that is, replacing $\lambda=-1+\gamma^2$ in \eqref{e:source}.
\end{Example}

\begin{Example}\label{ex:1ctd}
 Returning to Example \ref{ex:1}, we consider the (explicit) case of layers $w_*(x)=\tanh(x/\sqrt{2})$ connecting $w_\pm=\pm 1$ at $x=\pm\infty$. The eigenvalue problem  $w_{xx}+(1-3\tanh^2(x/\sqrt{2}))w=\lambda w$ can be converted into the first order system $u_x=A(x;\lambda)u$ with asymptotic matrices $A_\pm(\lambda)=\begin{pmatrix} 0 &1\\ \lambda+2 & 0\end{pmatrix}$. We have $i_\infty=1$ and stable and unstable subspaces are well defined outside of $\{\lambda\leq -2\}$. Solving the ODE explicitly, one finds the solution, substituting $\gamma=\sqrt{\lambda+2}$, 
 \[
 u_1^\mathrm{u}(x)= \begin{pmatrix} u_+(x)\\ u_+'(x) \end{pmatrix},\qquad
  u_2^\mathrm{s}(x)= \begin{pmatrix} u_+(-x)\\ -u_+'(-x) \end{pmatrix},
 \] 
 where
 \[
 u_+(x)= 
 (1 + \rme^{\sqrt{2}  x})^2
 \rme^{ 
 -\frac{
  \sqrt{2}\gamma (\sqrt{2}  - 3 \gamma + \sqrt{2}  \gamma^2)
  }{
  2 - 3 \sqrt{2}  \gamma + 2 \gamma^2
  }
   x 
}
  (2 - 3 \sqrt{2}  \gamma + 2 \gamma^2 + 
   4 \rme^{\sqrt{2}  x} (-2 + \gamma^2) + 
   \rme^{2 \sqrt{2}  x} (2 + 3 \sqrt{2}  \gamma + 2 \gamma^2)),
 \]
such that 
\[
 \iota(\lambda)=\begin{pmatrix} u_1^\mathrm{u}(0)&u_2^\mathrm{s}(0)\\
                  (u_1^\mathrm{u})'(0)&(u_2^\mathrm{s})'(0)
                \end{pmatrix}
                = \begin{pmatrix} -1+2\gamma^2 & -1+2\gamma^2 \\
                   -2 \gamma (-2 + \gamma^2)&2 \gamma (-2 + \gamma^2) 
                  \end{pmatrix},\qquad
\mathcal{E}(\lambda)=\mathrm{det}(\iota(\lambda))=-4 \gamma (-2 + \gamma^2) (-1 + 2 \gamma^2).
\]
Clearly, $\iota$ is analytic in $\gamma\in \C$ in this case, with zeros alias eigenvalues at $\gamma = 0,\pm \sqrt{2}, \pm 1/\sqrt{2}$. Only positive values of $\gamma$ correspond to eigenfunctions, negative values to resonance poles (exponentially growing solutions) and $\gamma=0$ to an embedded eigenvalue at the edge of the essential spectrum. Note that all roots of $\mathcal{E}$ are simple in this case. We see that $\iota$ is analytic on the Riemann surface defined by $\gamma$. 
We emphasize that our choice of $u_+(x)$ is by no means unique. One can clearly multiply $u_1^\mathrm{u}$ and $u_2^\mathrm{s}$ by non-vanishing analytic functions $\alpha_\pm(\lambda)$. In fact, canonical computations of the bases may well lead to choices where $\alpha_\pm(\lambda)$ have poles in the complex plane, which one then simply removes multiplying by suitable polynomials. A simple example of such a scaling is when one insists on a normalization $u_+(0)=1$, introducing a singularity $(1-2\gamma^2)^{-1}$ with two poles.  Less fortunate choices may introduce factors that exhibit additional branch points or other singularities, in the parametrization. An example for such a difficulty arises when attempting the common normalization $\mathcal{E}\to 1$ for $\lambda\to\infty$, which one could accomplish by normalizing $u_+(0)=(-1+2\gamma^2)/\gamma^{5/2}$, clearly introducing additional branch singularities. Another natural choice of normalization would be  $|u_1^\mathrm{u}(0)|=1$, which would, in addition to singularities, introduce terms involving $\bar{\gamma}$, destroying analyticity entirely.
\end{Example}
\begin{Example}[Lack of continuity]\label{ex:cpw}
 In function spaces, one readily concludes that invertibility is an open property in the spectral parameter, also under large classes of perturbations, which establishes upper semicontinuity of the spectrum under perturbations. This is, in general, not true for singularities of the pointwise resolvent as can be seen in the following example, borrowed from \cite{holzerscheel14},
 \begin{equation}\label{e:cpw0}
 u_t=-u_x+\eps v,\qquad v_t=v_x,
 \end{equation}
 which leads to the first order spatial spectral ODE
 \begin{equation}\label{e:cpw0s}
 u_x=-\lambda u+\eps v,\qquad v_x=\lambda v,
 \end{equation}
and globally analytic stable and unstable subspaces,
\[
E_+^\mathrm{s}(\lambda)=\begin{pmatrix}1\\0 \end{pmatrix},\qquad 
E_-^\mathrm{u}(\lambda)=\begin{pmatrix}\eps\\2\lambda \end{pmatrix},
\]
that intersect nontrivially at $\lambda=0$, $\mathcal{E}(\lambda)=2\lambda$. For $\eps=0$, however, the basis of $E_-^\mathrm{u}(\lambda)$ is degenerate at $\lambda=0$ so that a reparametrization is needed, for instance
\[
E_+^\mathrm{s}(\lambda)=\begin{pmatrix}1\\0 \end{pmatrix},\qquad 
E_-^\mathrm{u}(\lambda)=\begin{pmatrix}0\\1 \end{pmatrix}.
\]
As a result, the intersection is always trivial and $\mathcal{E}(\lambda)=1$. Put in the context of perturbation theory, the pointwise resolvent does not have a singularity for $\eps=0$, but upon arbitrarily small perturbations, such a singularity can be created.

The effect is of course also visible in the (explicit) solution to the equation, which for $\eps=0$ simply advects compactly supported initial conditions to the left ($u$-equation) and to the right ($v$-equation), which constitutes an effective super-exponential pointwise decay to zero. Coupling with $\eps\neq 0$ causes $u$ to converge to a constant, effectively integrating the initial mass in the $v$-equation. The effect appears also in less obvious examples, including for instance diffusion in \eqref{e:cpw0} or more general coupled amplitude equations \cite{faye17}.

We return to this example in \S\ref{s:5}, demonstrating how our algorithm correctly identifies the subtle dependence on the presence of a coupling term. 
\end{Example}

\begin{Example}[Branch poles vs branch points]\label{ex:bp}
In the trivial example $w_t=w_{xx}$, one finds $E_+^\mathrm{s}(\lambda)=(1,\sqrt{\lambda})^T$, $E_-^\mathrm{u}(\lambda)=(1,-\sqrt{\lambda})^T$, so $\mathcal{E}(\lambda)=2\sqrt{\lambda}$, which is \emph{both} not analytic at $\lambda=0$ due to a branch point in the eigenspaces and vanishes, so that $\iota^{-1}$ possesses a singularity of type $\sqrt{\lambda}^{-1}$. Passing to the Riemann surface by introducing $\gamma=\sqrt{\lambda}$, corresponding to considering $u_{tt}=u_{xx}$, one finds a simple pole at $\gamma=0$. 

Considering  $w_t=w_{xx}$ in $x>0$ with Robin boundary condition $n_1 w + n_2 w_x=0$ at $x=0$, one forms the Evans function from $E_\mathrm{bc}=(n_2,-n_1)^T$ and $E_+^\mathrm{s}(\lambda)=(1,\sqrt{\lambda})^T$ so that $\mathcal{E}(\lambda)=n_2\sqrt{\lambda}-n_1$, which still possesses a branch point singularity at $\lambda=0$, but does not vanish when $n_1\neq 0$. On the Riemann surface, we find a root $\gamma=n_1/n_2$, which corresponds to an eigenvalue when $n_1n_2>0$ and to a resonance otherwise.

\end{Example}
We refer to \cite{MR3100266} for many more examples and context.

% 
% example: iota for allen-cahn, no singularities on Riemann surface... see \begin{verbatim}lin_ac_layer.nb\end{verbatim} 
% 
% another example: where we should create an eigenvalue through the exponential decay! something with decay exp(-x), bc at x=0, cpw?

\subsection{Determinants and numerical methods}
We briefly comment on other numerical approaches related to this pointwise formulation with the aim of differentiating our approach from others in the literature. Finding spectral values, that is, points $\lambda$ where the inverse of $\iota(\lambda)$ is not analytic, can be reduced to taking a determinant of $\iota$ and finding roots of the resulting analytic function --- after first identifying branch points as a source of non-analyticity in the far field. For this, one needs to overcome several obstacles, starting with the computation of analytic bases in stable and unstable subspaces. One can track subspaces using differential forms, at the expense of a possibly high-dimensional system, or computing orthogonalized stable bases, at the expense of loosing analyticity; see for instance \cite{MR2221065} and references therein. Analyticity can be restored on the level of a determinant \cite{MR2253406,MR2676976}, thereby yielding efficient methods for computing subspaces and finding eigenvalues through winding number computations \cite{MR3413592}. In fact, from this point of view the pointwise nature of the computation can be relaxed to improve numerical stability, still exploiting a determinant formulation and computing winding numbers \cite{MR3157977}. There do not appear to be algorithms that do not involve a separate treatment of core and farfield, and most algorithms rely to some extent on determinants and winding number computations. In contrast, the approach that we present in the next section, treats core and farfield simultaneously and avoids determinants and winding numbers altogether, thus presenting a useful ad hoc tool for the initial study of stability problems.

\section{Inverse power methods for locally analytic operator pencils}\label{s:3}

Motivated by the previous derivation of nonlinear eigenvalue problems, we study families of matrices $\iota(\lambda)\in\C^{N\times N}$, in a domain $\lambda\in U\subset \C$, and wish to find values $\lambda_*$ such that the inverse $\iota(\lambda)^{-1}$ is not analytic at $\lambda=\lambda_*$. We assume that $\lambda$ is meromorphic on a Riemann surface, that is, $\iota(\varphi(\gamma))$ is meromorphic in $\gamma$, where $\varphi$ resolves potential branch points. We do not assume that $\varphi$ is a priori known.  There are many methods available that find poles of $\iota(\lambda)^{-1}$ in the case where $\iota$ is analytic; see in particular \cite{gutteltisseur} for a recent review. Many methods ultimately rely on particular polynomial interpolations of $\iota(\lambda)$ and subsequent root finding or linearization of the matrix pencil \cite{MR3144797}. Much of the suitability of a method depends on what is known about $\iota$, or, in other words, how it is actually computed. In our case, one usually starts computing $\iota$ at a fixed point $\lambda_0$, computing stable and unstable subspaces and choosing bases. The main difficulty now is to continue these bases to nearby values of $\lambda$ in an analytic fashion. A key obstacle is that a naive parametrization of the subspace as a graph over the reference subspace at $\lambda=\lambda_0$ may fail at isolated points, leading to singularities in $\iota$ induced by the parametrization, as exemplified in Example \ref{ex:1ctd} when normalizing $u_+(0)=1$. Alternatively, orthogonalizing bases for the parametrization destroys analyticity; see again Example \ref{ex:1ctd}.

% Most approaches eventually rely on computing determinants and ultimately exploit winding numbers in order to locate roots. 

Our approach relies on local power series from the graph parametrization, only, yet finds spectral values of $\iota$ even past the radius of convergence of the power series and potential singularities induced by the parametrization. The local power series, as we shall explain in the next chapter, is readily computable solving homological Sylvester equations.

To set up the analysis, we fix a reference value $\lambda_0$ with the goal of finding spectral values of $\iota(\lambda)$ closest to $\lambda_0$. We assume without loss of generality that $\lambda_0=0$ possibly redefining $\lambda$. We assume that the matrix function $\iota$ has a local expansion in a convergent power series with radius of convergence $R$,
\begin{equation}
 \iota(\lambda)=\sum_{k=0}^\infty \iota_k \lambda^k, \qquad |\lambda|< R.
\end{equation}
If $\iota_0$ is not invertible, $\lambda=0$ is already a spectral value and we therefore assume henceforth that $\iota_0$ is invertible. Consider then the infinite-matrix operator acting on infinite sequences $\underline{u}=(u_j)_{j=1,2,\ldots}$,
\begin{equation}\label{e:A}
 \mathcal{A}:\underline{u}\mapsto \mathcal{A}\underline{u},\qquad 
 (\mathcal{A}\underline{u})_j=
    \left\{\begin{array}{ll}                                       
                        -\iota_0^{-1}\left( \iota_1 u_1+\iota_2 u_2+\ldots \right) ,& j=1,\\
                        u_{j-1},& j>1.
           \end{array}\right., 
    \qquad \text{or }\mathcal{A}=
    \left(\begin{array}{cccc}
                -\iota_0^{-1} \iota_1 & -\iota_0^{-1} \iota_2 & -\iota_0^{-1} \iota_3 & \cdots\\
                1 & 0 & 0 & \cdots\\
                0 & 1 & 0 & \cdots\\
                0 & 0 & 1 & \cdots\\
                \vdots & \vdots & \vdots & \ddots
          \end{array}
    \right).
\end{equation}The form of $\mathcal{A}$ is motivated by the case where $\iota$ is a polynomial and $\mathcal{A}$ can act on finite sequences. The polynomial $\iota$ can then be thought of as the characteristic equation to a multi-term recursion, which in turn can be written as a first-order recursion in a higher-dimensional ambient space. Iterating $\mathcal{A}$ is, in this case, simply the inverse power method for this matrix representation. 

Eigenfunctions solve $\mathcal{A}\underline{u}=z\underline{u}$. Inspecting the components of this equation with $j>1$, we find $u_{j+1}=z^{-1} u_j$, so that $u_j=z^{-j}u_0$ for some vector $u_0\in\C^N$. Setting $\lambda=z^{-1}$, the first equation in  $\mathcal{A}\underline{u}=z\underline{u}$ gives
\[
 -\iota_0^{-1}\left( \iota_1\lambda +\iota_2\lambda^2+\ldots \right)u_0=z(\lambda u_0),
\]
which after multiplying by $\iota_0$ and rearranging gives
\[
 \iota(\lambda)u_0=0.
\]
In other words, we ``linearized'' the nonlinear matrix pencil, that is, spectral values $\lambda$ of the nonlinear pencil $\iota$ now correspond to spectral values $z=\lambda^{-1}$ of the (regular) eigenvalue problem for $\mathcal{A}$. 

To access regular spectral values, one now has access to traditional methods for eigenvalue problems. The idea we pursue here is to iteratively compute  $\mathcal{A}^k\underline{u}_0$ and expect that iterates grow with the spectral radius of $\mathcal{A}$, aligning with the eigenvector to the largest eigenvalue, for random initial vectors $\underline{u}_0$. Such convergence does depend on the nature of the spectrum of $\mathcal{A}$ and we will study three cases of interest in the subsequent three sections, characterized in terms of the spectral value of $\iota(\lambda)$ in the sense of Definition \ref{d:sv}:
\begin{enumerate}
 \item the singularity of   $\iota(\lambda)^{-1}$ closest  to $\lambda_0=0$ is a pole and lies within the radius of convergence $R$, \S\ref{s:3.1};
 \item the singularity of   $\iota(\lambda)^{-1}$ closest  to $\lambda_0=0$ is a pole and  lies within a ball where $\iota(\lambda)$ is meromorphic, \S\ref{s:3.2};
 \item the singularity of   $\iota(\lambda)^{-1}$ closest  to $\lambda_0=0$ is a branch point singularity, \S\ref{s:3.3}.
\end{enumerate}

\subsection{Isolated point spectrum}\label{s:3.1}
Clearly, $\mathcal{A}$ is a rank-1, hence compact perturbation of the right-shift operator, so that one can readily compute Fredholm properties in typical function spaces explicitly. Defining for instance $\ell^p_\rho$ for $\rho>0$ as the space of sequences such that $(u_j \rho^{-j})_j\in \ell^p$, we find 
\[
 \mathrm{spec}_{\mathrm{ess}, \ell^p_\rho}(\mathcal{A})=\{|z|\leq \rho^{-1}\}.
\]
On the other hand, the first row $\mathcal{A}_1:\ell^p_\rho\to \R$ is bounded only when $\rho<R$. Choosing $\rho$ arbitrarily close to $R$, we can thereby find eigenvalues of $\mathcal{A}$ within $\{|z| >R\}$ as point spectrum. Equivalently, any spectral value $\lambda$ of the operator pencil $\iota(\lambda)$ that lies within the radius of convergence of the power series can be found as an eigenvalue in the point spectrum of $\mathcal{A}$ in an appropriately chosen weighted space. In particular, if $\iota(\lambda)$ possesses a spectral value $\lambda$ with $|\lambda|<R$, the power method applied to $\mathcal{A}$ generically identifies the smallest eigenvalue of $\mathcal{A}$. 

\begin{Proposition}[Inverse Power Method --- point spectrum within radius of convergence]\label{p:pm1}
Assume that the nonlinear matrix pencil $\iota(\lambda)$ with radius of convergence $R>0$ possesses a unique smallest spectral value $\{|\lambda_0|<R\}$.
% , with $\frac{\rmd}{\rmd\lambda}\mathrm{det}(\iota(\lambda))|_{\lambda=\lambda_0}\neq 0$. 
In particular, $\iota(\lambda)^{-1}$ is analytic in  $|\lambda|< |\lambda_0|+\delta,\,\lambda\neq \lambda_0$, for some $\delta>0$. Then the associated inverse power iteration
 \[
  \underline{u}_{k+1}= \mathcal{A}\underline{u}_k,
 \]
 defined on $\ell^p_\rho$ with $1\leq p\leq \infty$ and $|\lambda_0|<\rho <R$ converges for initial vectors $\underline{u}_0$ in the complement $V$  of a strict subspace of $\ell^p_\rho$ to eigenvalue and eigenvector in the sense that
 \[
  \underline{u}_k/|\underline{u}_k| \to \underline{u}_*, \qquad \mathcal{A}\underline{u}_*=\lambda_0^{-1} \underline{u}_*,\qquad \iota(\lambda_0)(\underline{u}_*)_1=0.
  \]    
In particular, $V$ contains sequences $\underline{u}$ with $\underline{u}_j=0,j\geq 2$ and $\underline{u}_1\in V_0$, the complement of a strict subspace of $\C^N$. 
\end{Proposition}
\begin{Remark}\label{r:p1}
 \begin{enumerate}
%   \item Given the eigenvector, it is of course easy to reconstruct the eigenvalue $z_0=\lambda_0^{-1}$ of $\mathcal{A}$. 
  \item By the Analytic Fredholm Theorem, eigenvalues of $\mathcal{A}$ in $\{|z|< \rho^{-1}\} $ are isolated and of finite   algebraic multiplicity. Shifting $\lambda\mapsto \lambda-\lambda_\mathrm{s}$ by a small generic shift would therefore guarantee that the assumption of the proposition holds.
  \item Straightforward extensions of this result can establish that iteration of generic two-dimensional subspaces yield the eigenspace of $\mathcal{A}$ to the two smallest eigenvalues, showing as a consequence the convergence of a $QR$-type iteration scheme. 
  \item The rate of convergence can be readily obtained from the proof as the ratio between $\lambda_0$ and the next-smallest spectral value $\lambda_1$. We may compute for instance the sequence of approximate spectral values $\lambda_{0,k}$ via 
  \[
   \lambda_{0,k}^{-1}=\langle \underline{u}_{k+1},\underline{u}_k\rangle/\langle \underline{u}_{k},\underline{u}_k\rangle, 
  \]
  with, say, $\langle \underline{u},\underline{v}\rangle = (u_1,v_1)$, the standard complex scalar product in $\C^N$.  
  One finds from the proof below that $\underline{u}_k= \lambda_0^{-k}\underline{u}_*+\rmO(\lambda_1^{-k})$, so that 
  \begin{equation}\label{e:asyeig}
   \lambda_{0,k}^{-1}=\lambda_0^{-1}+\rmO((\lambda_1/\lambda_0)^{-k}).
  \end{equation}
  
  \end{enumerate}
\end{Remark}
\begin{Proof}
 By the analytic Fredholm theorem, we can decompose $X=\ell^p_\rho=X_0+X_1$ into $\mathcal{A}$-invariant subspaces so that $\mathcal{A}|_{X_0}=\lambda_0^{-1}\mathrm{id}+N$ with $N$ nilpotent, $X_0$ finite-dimensional, and  the spectral radius of $\mathcal{A}|_{X_1}$ is strictly less than $\lambda_0^{-1}$. Within $X_0$, we can analyze the iteration in Jordan Normal Form and find convergence of vectors to the eigenspace. The component in $X_1$ will decay exponentially due to the renormalization. 
 
It remains to show that choosing sequences with support on the first entry is sufficient to achieve growth. We therefore need to show that there exists a vector in the kernel of the adjoint $\mathcal{A}^*-z$ whose first component does not vanish. For any such vector $\underline{w}$, we quickly find, writing
$\iota^M(\lambda)=\sum_{\ell=0}^M\iota_\ell \lambda^\ell$,
\[
 w_j=\sum_{k=0}^{j-1} z^{j-1-k} \iota_k^T v_1 = z^{j-1} ((\iota^{j-1})^T(z^{-1})v_1,
\]
for some vector $v_1$. In order for $\underline{w}\in \ell^q_{\rho^{-1}}$, we need $w_j\rho^j\in\ell^q$, in particular $w_jz^{-j}\to 0$, so that in fact $\iota(z^{-1})v_1=0$, that is, $v_1$ belongs to the kernel of the adjoint. 
Clearly, $w_j=0$ for all $j$ if $v_1=0$, so that for a nontrivial element in the kernel $v_1\neq 0$ and therefore $w_1=\iota_0^Tv_1\neq 0$ using invertibility of $\iota_0$. This concludes the proof.  
\end{Proof}

\subsection{Extended point spectrum}\label{s:3.2}
We now turn to the case where $\iota(\lambda)$ does not have spectral values in $\{|\lambda|<R\}$. We assume however here that $\iota(\lambda)$ does have a meromorphic continuation in $\{|\lambda|<M\}$ and a spectral value in this disk. Note that, by uniqueness of the extension of $\iota$, the notion of spectral value in this larger disk is well defined, while the notion of eigenvalue for the associated operator $\mathcal{A}$ is not well defined since infinite sums do not converge when substituting a potential eigenvector to an eigenvalue with $|z|>R$ into the expression for the first component $(\mathcal{A}\underline{u})_1$.

\begin{Proposition}[Inverse Power Method --- point spectrum within meromorphic domain]\label{p:pm2}
 Assume that the nonlinear matrix pencil $\iota(\lambda)$ is meromorphic in $|\lambda|<M$ possesses a unique smallest spectral value with $\{|\lambda_0|<M\}$, that is, $\iota(\lambda)^{-1}$ is analytic in  $|\lambda|< M,\,\lambda\neq \lambda_0$. Then, for any $K\geq 1$,  the associated inverse power iteration 
 \[
  \underline{u}_{k+1}= \mathcal{A}\underline{u}_k,%/|\mathcal{A}\underline{u}_k|_{\ell^p_\rho},
 \]
 with compactly initial data, $(\underline{u}_0)_j=0$ for all $j>K$, converges for all initial vectors $(\underline{u}_0)_{1\leq j\leq K}\in \C^K$ except for a finite-codimension subspace, locally uniformly. More precisely, for any $K_0$,
 the restriction to the first $K_0$ components $R_{K_0}\underline{u}=(u_1,\ldots,u_{K_0})$ converges to the restriction of a formal eigenvector,
 \[
  R_{K_0}\underline{u}_k/|R_{K_0}\underline{u}_k|
  \ \to R_{K_0}\underline{u}_*,
%   \qquad \mathcal{A}\underline{u}_*=\lambda_0^{-1} \underline{u}_*,\qquad \iota(\lambda_0)(\underline{u}_*)_1=0.
  \]
and 
\[
 R_{K_0}(\mathcal{A}\underline{u}_k-\lambda_0^{-1} \underline{u}_k )\to 0, \quad \text{ for } k\to\infty.
\]
\end{Proposition}
\begin{Remark}\label{r:p2}
 \begin{enumerate}
  \item Similar to the comments in Remark \ref{r:p1}, one can generalize to multiple leading eigenvalues using iteration of subspaces with appropriate orthogonalization strategies. 
  \item Convergence is again exponential, with rate given by the ratio between $\lambda_0$ and the next-smallest spectral value $\lambda_1$ as in \eqref{e:asyeig}.
 \end{enumerate}
\end{Remark}
To prepare for the proof, we introduce a pointwise description of iterates.  We wish to obtain a pointwise representation of $\mathcal{A}^k$, that is, for the matrix entries $((\mathcal{A}^k \delta_{jm})_\ell= ((\mathcal{A}^k)_{\ell m}$ for fixed $\ell$ and $m$. We wish to use Dunford's resolvent identity and  start with an expression for the resolvent $(z-\mathcal{A})^{-1}$. We therefore fix $m$ arbitrary and solve
 \[
  \left((z-\mathcal{A})\underline{u}\right)_m=f, \qquad   \left((z-\mathcal{A})\underline{u}\right)_j=0, \ j\neq m, 
 \]
explicitly.
We find, solving the equation for all $j>1$,
\begin{equation}\label{e:pr1}
 u_j=z^{-j}u_0,\ j<m,\qquad u_j=z^{-j}u_0 + z^{m-j-1}f,\ j\geq m.
\end{equation}
Inserting into the equation for $m=1$ gives 
\begin{align*}
0&= -\iota_0^{-1}\left(\iota_1 z^{-1}+\iota_2 z^{-2}+\ldots\right)u_0-u_0
 -\iota_0^{-1}\left(\iota_m z^{-1}+\iota_{m+1} z^{-2}+\ldots\right) f\\
 &=\iota(\lambda)u_0 - \lambda^{1-m}\left(\iota(\lambda)-\iota^{m-1}(\lambda)\right)f,
\end{align*}
where $\iota^{p}(\lambda)=\iota_0+\ldots+\iota_p\lambda^p$ is the Taylor jet up to order $p$. Solving this matrix equation with matrix entries in the field of meromorphic functions for $u_0$ gives 
\begin{equation}\label{e:pr2}
 u_0=\lambda^{1-m}\iota(\lambda)^{-1}\left(\iota(\lambda)-\iota^{m-1}(\lambda)\right)f,
\end{equation}
which together with \eqref{e:pr1} defines the pointwise resolvent $u_j=\mathcal{R}(z;\mathcal{A})_{jm}f$ when the right-hand side is supported in the $m$'th component. We write $\mathcal{R}(z;\mathcal{A})$ for the infinite matrix $1\leq j,m<\infty$. 

From the form of \eqref{e:pr1}--\eqref{e:pr2}, we obtain the following lemma.
\begin{Lemma}
 The pointwise resolvent  $\left((z-\mathcal{A})^{-1}\right)_{jk}$ possesses an analytic extension into connected component of the region $\{z=1/\lambda\}$ where $\iota(\lambda)$ is meromorphic and $\iota(\lambda)^{-1}$ is analytic. Moreover, if $\iota(\lambda)^{-1}$ has a pole at $\lambda_0$, then the components  $\left((z-\mathcal{A})^{-1}\right)_{j1}$ of the pointwise resolvent have a singularity at $z=z_0$. 
\end{Lemma}
\begin{Proof}
 We only need to show that the pointwise resolvent cannot be analytic when $\iota(\lambda)^{-1}$ is not analytic. This follows by setting $m=1$ in \eqref{e:pr2} so that, with \eqref{e:pr1},
 \[
  u_j=\lambda^{j}\left(\mathrm{id}-\iota(\lambda)^{-1}\iota_0\right)f.
 \]
Here, the term $\lambda^{j}f$ is analytic, and the term $\lambda^{j}\iota(\lambda)^{-1}\iota_0f$ has a singularity since $\iota_0\lambda^{j}$ is invertible. 
\end{Proof}

From the form of \eqref{e:pr1}--\eqref{e:pr2}, it is clear that the pointwise resolvent possesses an analytic extension into the region where $\iota(\lambda)^{-1}$ is analytic and $\iota(\lambda)$ is meromorphic. 
% We say that $\lambda$ lies in the \emph{extended point spectrum} of $\iota(\lambda)$ when the pointwise resolvent is not analytic at  $z=1/\lambda$. 

\begin{Proof}[ of Proposition \ref{p:pm2}.]
Choosing a contour $\Gamma=\{|z|=R\}$ with $R$ large, oriented counter-clockwise, one obtains from Dunford's calculus  that 
\[
\underline{u}^k:= \mathcal{A}^k \underline{f}=\frac{1}{2\pi\rmi}\int_\Gamma z^k (z-\mathcal{A})^{-1}\underline{f} \rmd z.
\]
For $\underline{f}$ compactly supported, and evaluating both sides in a compact region $j\leq J$, we may deform the contour $\Gamma$ in the region where the pointwise resolvent $((z-\mathcal{A})^{-1})_{jk}$ is analytic, that is, within the region where it is meromorphic but outside of the extended point spectrum. We choose to deform the contour into $\tilde{\Gamma}=\Gamma_0\cup\Gamma_1$, where $\Gamma_1=\{|z|=R_2<|\lambda_0|^{-1}\}$ and $\Gamma_0=\{\lambda_0^{-1}+z|\,|z|=\eps\} $ for some sufficiently small $\eps>0$. For the contribution from $\Gamma_1$, one readily finds componentwise decay $|\underline{u}^k_j|\leq C R_2^k$. The contribution from $\Gamma_0$ can be evaluated computing residuals after expanding the pointwise resolvent in a Laurent series, which gives a contribution $\sum_{\ell=0}^{\ell_0} Q_j k^k \lambda_0^k$. From this splitting, the claim follows readily, in complete analogy to the finite-dimensional convergence of the power method.
\end{Proof}
\begin{Remark}[Zeros of meromorphic functions]\label{r:z}
 The strategy employed here can of course be most easily tested as an algorithm to find roots of meromorphic functions $f(\lambda)$ in the plane $z\in\C$. More precisely, our algorithm finds the zero $\lambda_*$ of $f(\lambda)$ closest to a fixed reference point $\lambda_0$ using only the Taylor expansion of $f$ at $\lambda_0$. One simply iterates 
 \[
 u_k=\frac{-1}{f(\lambda_0)}\left( f'(\lambda_0)u_{k-1}+ \frac{1}{2}f''(\lambda_0)u_{k-2}+ \frac{1}{6}f'''(\lambda_0)u_{k-3}+\ldots   \right), \qquad u_0=1,\  u_j=0\text{ for } j<0,
 \]
 and obtains $\lambda_*-\lambda_0=\lim_{k\to\infty} u_k/u_{k+1}$.  
 Our result here states that this iterative algorithm identifies zeros past the radius of convergence of the local power series. Of course, this approach is useful only when access to Taylor series coefficients is preferred to simple evaluation of a function. 
\end{Remark}

\subsection{Branch points}\label{s:3.3}
A third typical possibility appears when the largest singularity of $(z-\mathcal{A})^{-1}$ is a branch point singularity. We say that $\iota$ has a branch pole of order $p$ for some $p\in\N$ at $\lambda_0$ if $\iota(\lambda_0+\gamma^q)^{-1}$ is componentwise meromorphic in $\gamma$ near $\gamma=0$ with a simple pole at $\gamma=0$ for $p=q$, but is not meromorphic for $1\leq q<p$. We focus here on the case $p=2$. 
  
For any $\lambda_0\neq 0$, let $S_\theta(\lambda_0)$ be the sector $\{\lambda\,|\,\mathrm{arg}((\lambda-\lambda_0)/\lambda_0)<\theta\}$ and $B_r=\{\lambda|\,|\lambda|<R\}$. 

\begin{Proposition}[Inverse Power Method --- branch points within meromorphic domain]\label{p:a3}
Given $\lambda_0\neq 0$, $|\lambda_0|=M$, $\delta>0$, and $\theta<\pi/2$, define 
$\Omega =B_{M+\delta}\setminus \overline{S_\theta(\lambda_0)})$. 
Assume that the nonlinear matrix pencil $\iota(\lambda)$ is pointwise meromorphic  in $\Omega$ and has a branch pole of order 2 at $\lambda_0$.

% pointwise continuous in $\overline{\Omega}$, and $\mathrm{det}(\iota(\lambda)) \neq 0$ in $\overline{\Omega}$. 

Then the associated inverse power iteration 
 \[
  \underline{u}^{k+1}= \mathcal{A}\underline{u}^k,
 \]
 with compactly supported  initial data, $(\underline{u}_0)_j=0,j>K$ asymptotically exhibits pointwise exponential growth with rate $1/\lambda_0$ with an algebraic correction,
 \[
   \underline{u}^k_j = \lambda_0^{-k} k^{-1/2} P_j  \underline{u}^0 \left(1+\rmo_1(k^{-1})\right).
 \]
 for some non-vanishing linear map $P_j$ defined on compactly supported sequences. 
%  
%  There are $z_0$ and $0<R_1<|z_0|$ such that $(z-\mathcal{A})^{-1}$ is componentwise analytic in $\{|z|\geq R_1\}\setminus S_{z_0}$
%  
%  , for some sector $S_{z_0}=\{z| |\mathrm{arg}\,(z/z_0)-\pi|<\delta$ for some $0<\delta
\end{Proposition}
\begin{Remark}\label{r:p3}
  \begin{enumerate}
   \item For higher-order branch points with Riemann surface covering $\lambda=\lambda_0+\gamma^p$, one finds in an equivalent fashion asymptotics with growth $\lambda_0^{-k}k^{1-1/p}$. 
   \item Another case of interest arises in $x$-dependent problems when $\iota$ possesses a branch point singularity but $\iota^{-1}$ is continuous. In this case, for $p=2$, one finds pointwise rates $\lambda_0^{-k}k^{-3/2}$ in analogy to the pointwise decay for the heat equation on the half line with Dirichlet boundary condition. 
   \item From the asymptotics for $\underline{u}^k$ with rate $\lambda_0^{-k}k^{-\alpha}$, one readily derives asymptotics of $\lambda_{0,k}$ as in Remark \ref{r:p1} (iii),
   \begin{equation}\label{e:bplamasy}
      \lambda_{0,k}\sim\lambda_0+\frac{\alpha\lambda_0}{k}.
   \end{equation}
   In particular, predictions for the branch point converge algebraically, with rate $k^{-1}$, regardless of the order of the branch point and $\alpha$, but with a prefactor $\lambda_0$ which is small for good initial guesses, suggesting effective shift strategies. Iterating a finite number $K$ of iterates to find a new initial guess $\lambda^K_0$ and restarting with the new initial guess $\lambda_0^K$, one finds exponential convergence in $k$. We demonstrate this strategy in \S\ref{s:5}. 
  \end{enumerate}
\end{Remark}

\begin{Proof}
The inverse power operator $\mathcal{A}$ associated with $\iota$ is invertible in $z\in \overline{\Omega'}$, where $\Omega'=1/\Omega$ contains all inverses $\lambda^{-1}$ of elements in $\Omega$. 
We can therefore write, in a pointwise sense,
\[
\underline{u}^k= \mathcal{A}^k \underline{f}=\frac{1}{2\pi\rmi}\int_\Gamma z^k (z-\mathcal{A})^{-1}\underline{f} \rmd z, 
\] 
for $\Gamma=\partial\Omega'$. Here, we use that the singularity of $\iota(\lambda)$ at $\lambda_0$ due to the simple pole in $\gamma$ is integrable, $\rmO(\lambda^{-1/2})$, leading to an integrable singularity of $(z-\mathcal{A})^{-1}$ on $\Gamma$. 
In the following, we assume for simplicity that $\lambda_0=1$, the general case can be easily obtained from there by scaling and complex rotation.
Expanding the pointwise resolvent of $\mathcal{A}$ near $z_*=1/\lambda_0=1$, we write $(z-\mathcal{A})^{-1}=(z-1)^{-1/2}\mathcal{B}_0 +\rmO(1)$, which gives
\[
\underline{u}^k=\frac{1}{2\pi\rmi}\int_\Gamma z^k \left((z-1)^{-1/2}\mathcal{B}_0 +\rmO(1)\right)\underline{f} \rmd z, 
\]
Ignoring contributions from $\Gamma$ where $|z|<1-\delta$ for some $\delta>0$, we parameterize $\Gamma=\Gamma_+\cup\overline{\Gamma_+}$, with $\Gamma_+=\{z=1-\rme^{\rmi\theta\tau},0\leq \tau\leq \delta$, and find
\begin{align*}
 \underline{u}^k\sim &\frac{\rme^{\rmi\theta}}{2\pi\rmi}\int_0^\delta (1-\rme^{\rmi\theta}\tau)^k \left((-\rme^{\rmi\theta}\tau)^{-1/2}\mathcal{B}_0 +\rmO(1)\right)\underline{f} \rmd \tau -\frac{\rme^{-\rmi\theta}}{2\pi\rmi}\int_0^\delta (1-\rme^{-\rmi\theta}\tau)^k \left((-\rme^{-\rmi\theta}\tau)^{-1/2}\mathcal{B}_0 +\rmO(1)\right)\underline{f} \rmd \tau \\
 =& -\frac{1}{\pi}\int_0^\delta (1-\tau)^k(\tau^{-1/2}\mathcal{B}_0+\rmO(1))\underline{f}\rmd\tau=k^{-1/2}P\underline{f}\left(1+\rmo_1(k^{-1}\right).
\end{align*}
% 
% 
%  Invertibility of $\iota$ implies pointwise invertibility of $\mathcal{A}$ such that we find the pointwise estimate 
% \[
% | \underline{u}^k_j|\leq C\int_\Gamma |z^k|\rmd z.
% \]
% Clearly, as $k\to\infty$, the integral on parts of $\Gamma$ away from the inverse of the sector $S_\theta(\lambda_0)$ decay exponentially in $k$ so that it suffices to consider the integral along $\partial S_\theta(\lambda_0)$. Here, $|z|=1-\mu$, where $\mu$ parameterizes the edge of the sector and we readily find 
% \[
% | \underline{u}^k_j|\leq C_1\int_0^\delta |(1-\mu)^k|\rmd \mu\leq C_2/k.
% \]
\end{Proof}
% \begin{Remark}
% Assuming more regularity and precise conditions on the behavior of $\iota$ near the branch point, one can obtain asymptotic expansions, commonly starting with a term $k^{-3/2}$, again for almost all initial conditions supported on $j=1$. For this, one parameterizes $\Gamma$ and expands the (pointwise) resolvent near the pole, finding that even-order terms in the expansion cancels. 
% 
% One can similarly treat the case where $\iota(\lambda_0)$ is not invertible, leading to weaker algebraic decay. 
% \end{Remark}

\section{Implementation of algorithms}\label{s:4}

Practically,  we wish to start with an ``explicit'' matrix-valued family $A(x;\lambda)$ and asymptotic matrices $A_\pm(\lambda)$ as in \eqref{e:twlin}, all polynomial in $\lambda$. In order to apply the inverse power method as described above, we need to 
\begin{enumerate}
 \item find a basis for $E^\mathrm{u}_-(\lambda_0)$ and for $E^\mathrm{s}_+(\lambda_0)$;
 \item compute Taylor expansions for  $E^\mathrm{u}_-(\lambda)$ and for $E^\mathrm{s}_+(\lambda)$ at $\lambda=\lambda_0$;
 \item assemble the map $\iota(\lambda)$ represented by a power series and implement the inverse power iteration. 
\end{enumerate}
We describe these somewhat practical issues in the next three sections. 

\subsection{Finding invariant subspaces and computing Taylor jets}\label{s:inv}

We describe how to obtain invariant subspaces, expand in $\lambda$, and continue using Newton's method. 
\paragraph{Schur decomposition.} Typical starting point for spectral computations is the region where stable and unstable subspaces actually correspond to the $k$ most unstable and $N-k$ most stable eigenvalues, respectively. Of course, we are particularly interested in situations where this splitting is no longer valid at the relevant eigenvalue $\lambda$, but subspaces at these values are the analytic continuation from values where the splitting is valid. We use a Schur decomposition sorting by real parts of eigenvalues to find an orthonormal basis and an orthonormal complement to $E^\mathrm{u}_-$ and $E^\mathrm{s}_+$ from the matrices $A_\pm(\lambda_0)$, all arranged in orthonormal matrices $U^\mathrm{s/u}_\pm$.

\paragraph{Taylor jets.} Computing Taylor jets for subspaces is a special case of computing Taylor expansions for invariant manifolds, which one readily sees by appending the trivial equation $\lambda'=0$. We outline the relevant steps, here.  We first shift the polynomial pencil evaluating derivatives at $\lambda_0$ and then conjugate with $U^\mathrm{s/u}$ so that $(U^\mathrm{s/u})^TA_\pm(\lambda+\lambda_0) U^\mathrm{s/u}$ possesses the trivial invariant subspace spanned by the first $k$ or $N-k$ coordinate vectors at $\lambda=0$, respectively. In the following, we therefore outline how to compute expansions near $\lambda=0$ for a polynomial pencil of degree $p$ with block form corresponding to the decomposition $\C^N=E_0\oplus E_1$ into canonical eigenspaces,
\[
 A(\lambda)=\begin{pmatrix} 
              A_{00}(\lambda) & A_{01}(\lambda)\\ A_{10}(\lambda) & A_{11}(\lambda)
            \end{pmatrix},
\qquad A_{10}(0)=0, \quad A_{00}\ k\times k-\text{matrix}, \ A_{11}\ (N-k)\times (N-k)-\text{matrix}. \
\]
We write the invariant subspace as a graph of $H(\lambda):E_0\to E_1$, $H(0)=0$, giving the column representation $E^\mathrm{s/u}\sim U^\mathrm{s/u}(F_0+H(\lambda)F_0) (U^\mathrm{s/u})^T$, where the $N\times k$-matrix $F_0$ forms the canonical basis in $E_0$. Invariance of $\text{graph}(H)$, that is,
\[
 A(\lambda) \left\{ \begin{pmatrix} F_0\\H(\lambda)F_0\end{pmatrix},\ F_0\in E_0 \right\}=  \left\{ \begin{pmatrix} F_1\\ H(\lambda)F_1\end{pmatrix},\ F_1\in E_0 \right\},
\]
is equivalent to requiring that for each $F_0\in E_0$, there exists $F_1\in E_0$ so that
\[
 A(\lambda) \begin{pmatrix} F_0\\H(\lambda)F_0\end{pmatrix}=  \begin{pmatrix} F_1\\H(\lambda)F_1\end{pmatrix}.
\]
This gives the matrix identity
\begin{equation}\label{e:hom}
 A_{10}(\lambda)+A_{11}(\lambda)H(\lambda)=H(\lambda)A_{00}+H(\lambda)A_{01}(\lambda)H(\lambda).
\end{equation}
Expanding $H$ and the $A_{jk}$ in $\lambda$ via 
\[
A_{jk}(\lambda)=\sum_{\ell=0}^pA_{jk}^\ell \lambda^\ell,\qquad H(\lambda)=\sum_{\ell=0}^\infty H^\ell \lambda^\ell,
\]
we find that $A_{10}^0=0,\ H^0=0$, and, at order $\ell$, 
\begin{equation}
A_{11}^0 H^\ell - H^\ell  A^0_{00} = R^\ell, \qquad R^\ell= \sum_{j=1}^{\ell-1}\left( H^jA_{00}^{\ell-j} -  A_{11}^{\ell-j} H^j \right) -A_{10}^\ell + \sum_{\substack{i+j+k=\ell\\0\leq j\leq p\\1\leq i,k\leq \ell-1}}H^i A_{01}^jH^k.
\end{equation}
At each order $\ell=1,2,\ldots$, this equation can be solved for $H^\ell$ by solving a linear Sylvester equation for $H^\ell$, with linear operator explicit on the left-hand side. The Sylvester equation can be solved effectively putting $A_{00}$ and $A_{11}$ into upper triangular form using Schur decomposition. For finite (low) order $p$, the right-hand side requires $\rmO(\ell)$ matrix multiplications so that overall effort is quadratic in the maximal order $\ell$.

\paragraph{Newton's method and continuation.}
We note that the formulation here also lends itself to direct Newton and continuation approaches, which we shall exploit when restarting the inverse power iteration. An approximate invariant subspace solves \eqref{e:hom} for some $\lambda_*$ with a small residual. Using Newton's method, solving again a Sylvester equation at each step, we can find a nearby actual invariant subspace. We can also implement continuation in $\lambda$, choosing for instance a generic complex path between two spectral parameter values $\lambda_0$ and $\lambda_1$ of the form 
\[
\lambda(\tau)=\lambda_0+\tau(\lambda_1-\lambda_0)+ \rmi \rho (\lambda_1-\lambda_0)\tau(1-\tau), \quad \rho\in [-1,1] \text{ fixed}.
\]
For a generic choice of $\rho$, the path would avoid isolated poles of $H$ or branch point singularities of the subspace so that arclength continuation would successfully find the desired invariant subspace at $\lambda_1$, even if that subspace is not actually the unstable subspace. 

% * references: Beyn?

\subsection{Assembling $\iota$}
We illustrate how to assemble $\iota$ in the simple case of a discretization based on the second order trapezoidal rule. Let $(u_j)_{j=1\ldots n+1}$ be the values at grid points $x_j$ and $u_\mathrm{bc}=(u_\mathrm{u},u_\mathrm{s})\in\C^k\times \C^{N-k}$ a vector parameterizing boundary conditions. The differential equation is then encoded in the $Nn\times N(n+2)$-matrix corresponding to  $\frac{1}{h} (u_{j+1}-u_j)=\frac{1}{2}(A(x_{j+1};\lambda)+ A(x_j;\lambda))$, with zero columns at the end corresponding to $u_\mathrm{bc}=(u^\mathrm{u},u^\mathrm{s})$. We add $2N$ rows corresponding to $u_1=U^\mathrm{u}(\lambda)u^\mathrm{u}$ and $u_1=U^\mathrm{s}(\lambda)u^\mathrm{s}$, where $U^\mathrm{s/u}(\lambda)$ are bases for $E^\mathrm{u/s}_\pm(\lambda)$. The resulting $N(n+2)\times N(n+2)$ square matrix is the desired nonlinear matrix family $\iota(\lambda)$. It is sparse at any order $\lambda$ with entries in $N\times 2N$ blocks along the diagonal at orders $\ell\leq p$ and with nonzero entries only in the bottom right $2N\times N$-corner for orders $\ell>p$.

For constant coefficients, the differential equation can of course be ignored and $\iota$ is simply given by the $N\times N$-matrix $(U^\mathrm{u}(\lambda)|U^\mathrm{s}(\lambda))$.

We implemented the family $\iota(\lambda)=\iota^0+\iota^1\lambda+\ldots $ as a sparse matrix $\iota=(\iota^0|\iota^1|\iota^2|\ldots|\iota^M)$ allowing easy extraction of orders of iota for the inverse power iteration.

\subsection{Implementing the inverse power method}

We initiate the inverse power iteration iterating $\mathcal{A}$ in \eqref{e:A} with a random complex starting $N$-vector $u_1$. Note that the method involves shifting only, in all but the first component. In the first component, we apply the pencil expansion terms $\iota_\ell$ and solve a linear equation with matrix $\iota_0$.
Having precomputed expansions up to an order $M$, we can then perform $M$ iterates exactly. Predictions for the eigenvalue are obtained from the first component $\lambda_\mathrm{p}=\langle u_1, u_1\rangle/\langle u_1, u_2\rangle$. Stopping criteria are formulated in terms of tolerances for the change in $\lambda_\mathrm{p}$ and the first components $\|\lambda_\mathrm{p}u_2-u_1\|$. After $M$ iterations or when initial tolerances are met, we restart the pencil iteration: we shift the symbol $\iota$ to the new predicted value $\lambda_\mathrm{p}$, shifting polynomials explicitly and recomputing eigenspaces using either continuation or a Newton method with predictor from Taylor expansion, as described in \S\ref{s:inv}. For these subsequent iterations, we use a lower truncation order of the pencil $M_\mathrm{fine}\ll M$ with frequent restarts until a fine tolerance is met. Shifts using step sizes roughly $\tau(\lambda_\mathrm{p}-\lambda_\mathrm{old})$ with $\tau\sim 0.8\ldots 0.95$ turn out to be most robust avoiding both the problem of non-invertibility of $\iota_0$ at the sought-after eigenvalue and problems of continuing and computing eigenspaces at branch points. 
% We also implemented standard $QR$-iteration to compute multiple spectral values simultaneously.

Since convergence near branch points is slow, algebraic, we also implemented a Newton method to find the exact location of branch points for constant coefficient problems. Branch points solve the system
\begin{align*}
\begin{array}{rrr} A(\lambda)u-\nu u=0,& \qquad\qquad \qquad  &
\langle e_0,u\rangle -1=0,\\
A(\lambda) v-\nu v-u=0,&\qquad \qquad\qquad  &
\langle e_0,v\rangle =0.\,
\end{array}
\end{align*}
where $e_0$ is an approximate element of the kernel of $A(\lambda)-\nu$ and the scalar products are understood as Hermitian (complex valued) forms. The inverse power iteration provides good initial guesses for $\lambda$. We find an initial guess for $u$  by computing the intersection of $E^\mathrm{s}_+$ and $E^\mathrm{u}_-$ at the initial guess and computing eigenvalues $\nu$ and eigenvectors $u$ for $A(\lambda)$ restricted to this intersection.

\section{Numerical examples}\label{s:5}

We demonstrate convergence and effectiveness of the algorithms in several examples. 

\paragraph{Pointwise growth modes --- constant coefficients and branch points of the dispersion relation.}

\begin{figure}%[h!]
\includegraphics[width=0.33\textwidth]{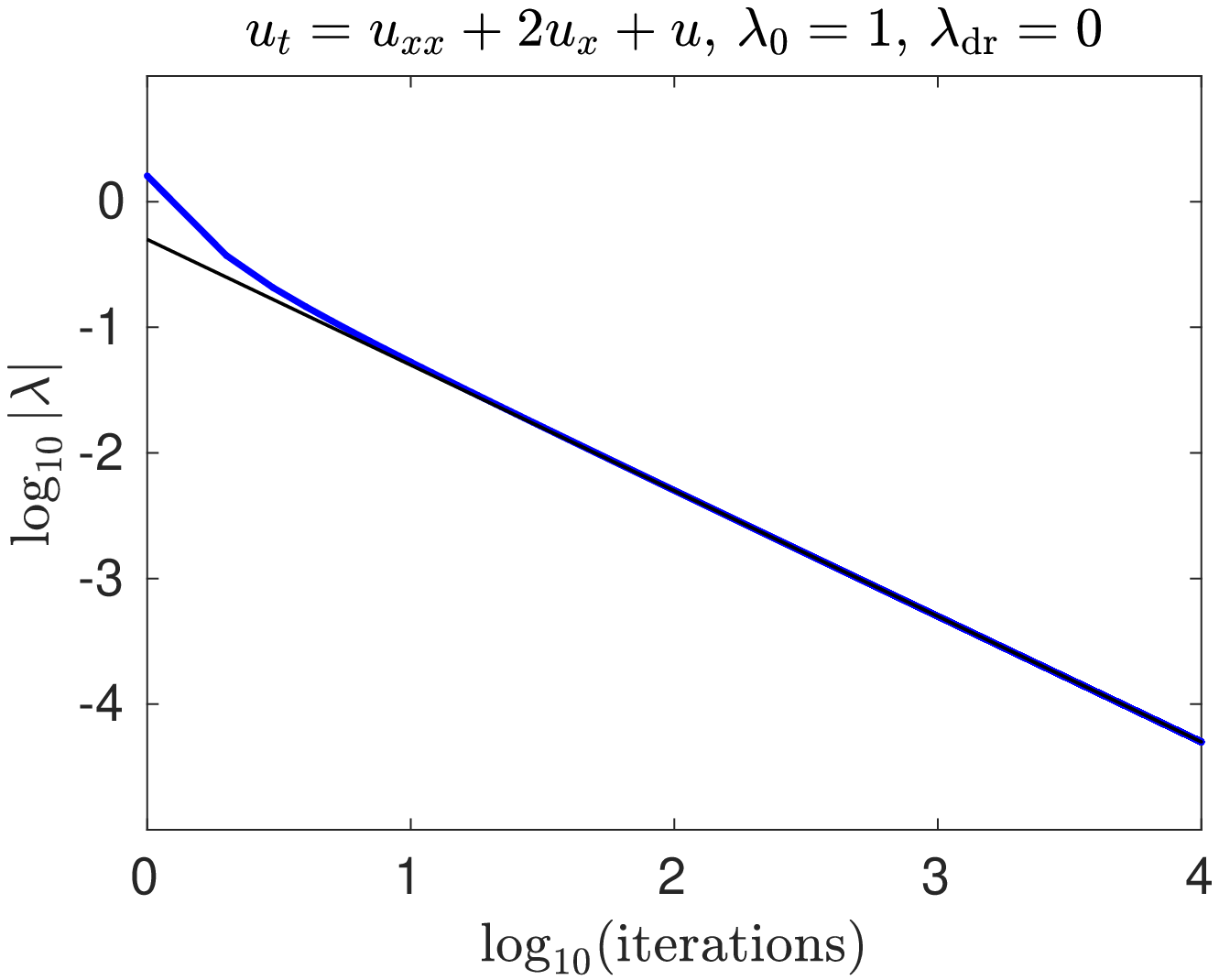}%
\includegraphics[width=0.33\textwidth]{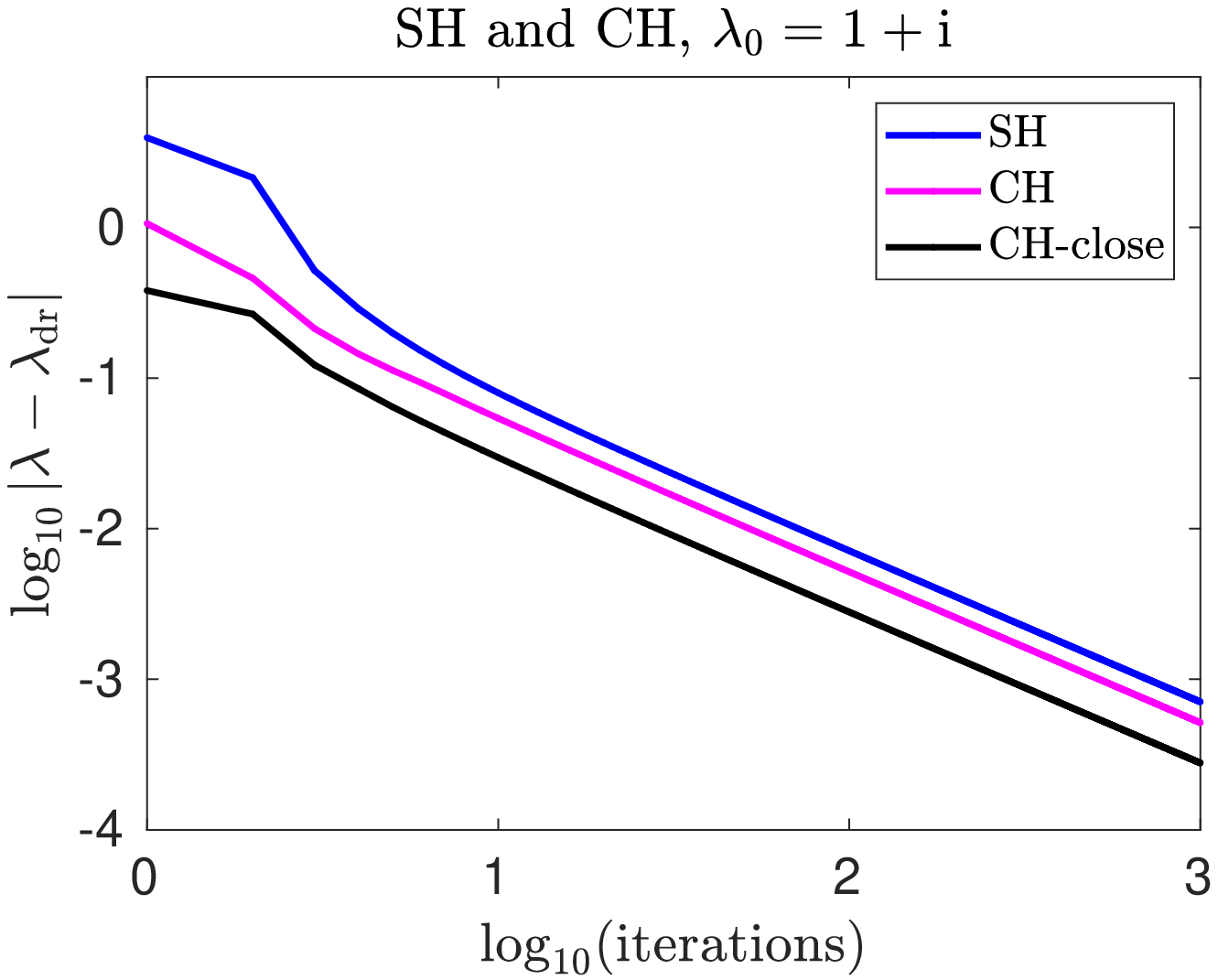}%
\includegraphics[width=0.33\textwidth]{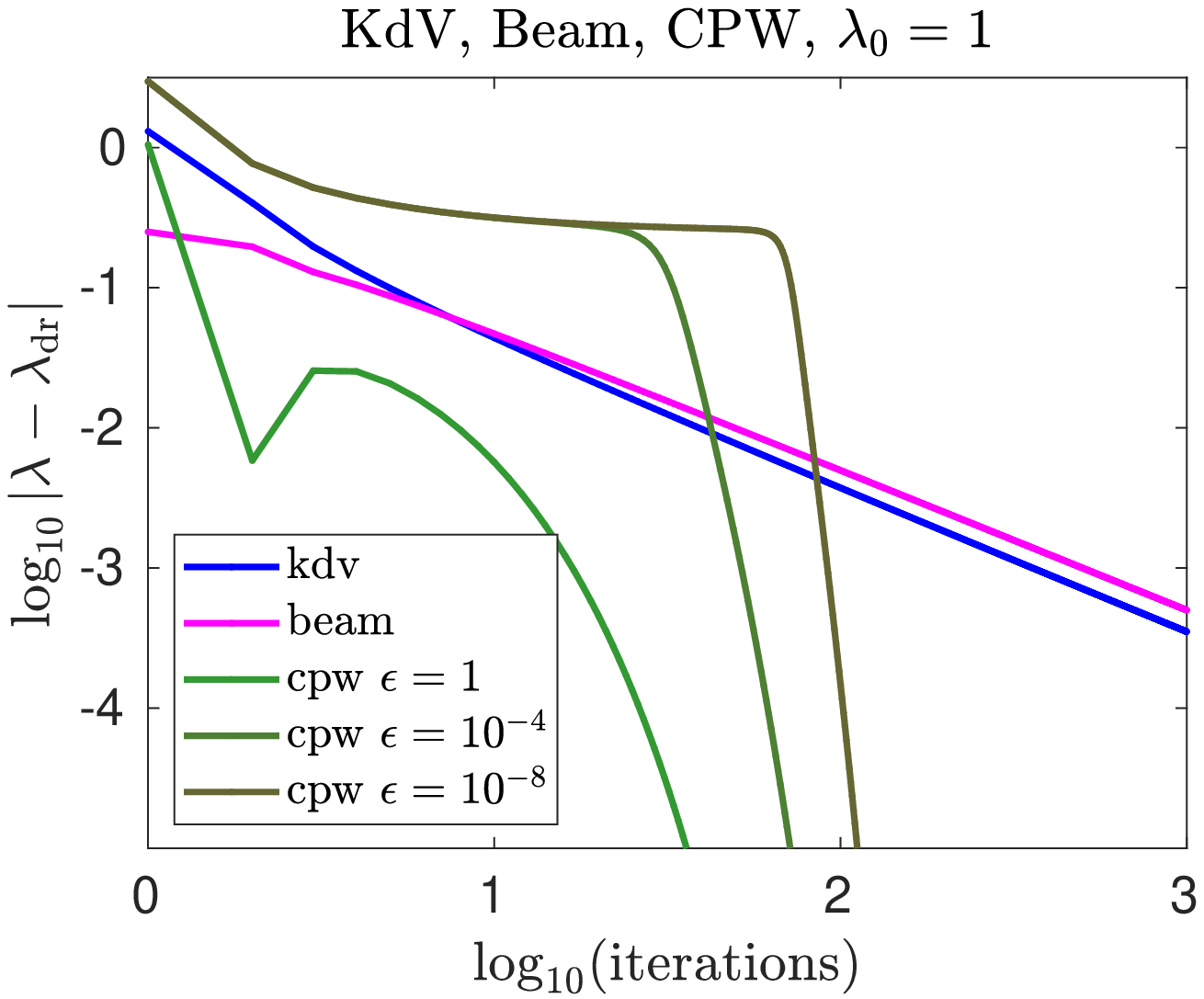}%
\caption{\emph{Left:} Convergence to $\lambda_\mathrm{dr}=0$ in convection-diffusion \eqref{e:cd} with starting value $\lambda_\mathrm{0}=1$ and linear fit with slope $-1$ corresponding to an algebraic convergence rate $k^{-1}$. \emph{Center:} Convergence to $\lambda_\mathrm{dr}=0$ in the Swift-Hohenberg (SH) equation \eqref{e:sh} and to $\lambda_\mathrm{dr}=\rmi\omega_\mathrm{dr}$ in the Cahn-Hilliard (CH) equation \eqref{e:ch} with starting values $1+\rmi$ and $0.5+\rmi$ (CH only). \emph{Right:}  Algebraic convergence to $\lambda_\mathrm{dr}=0$ for multiple double roots in KdV \eqref{e:kdv} and beam equation \eqref{e:beam}, as well as exponential convergence in the coupled transport equation (CPW); see text for details.  }\label{f:1}
\end{figure}

\begin{figure}
\includegraphics[width=0.33\textwidth]{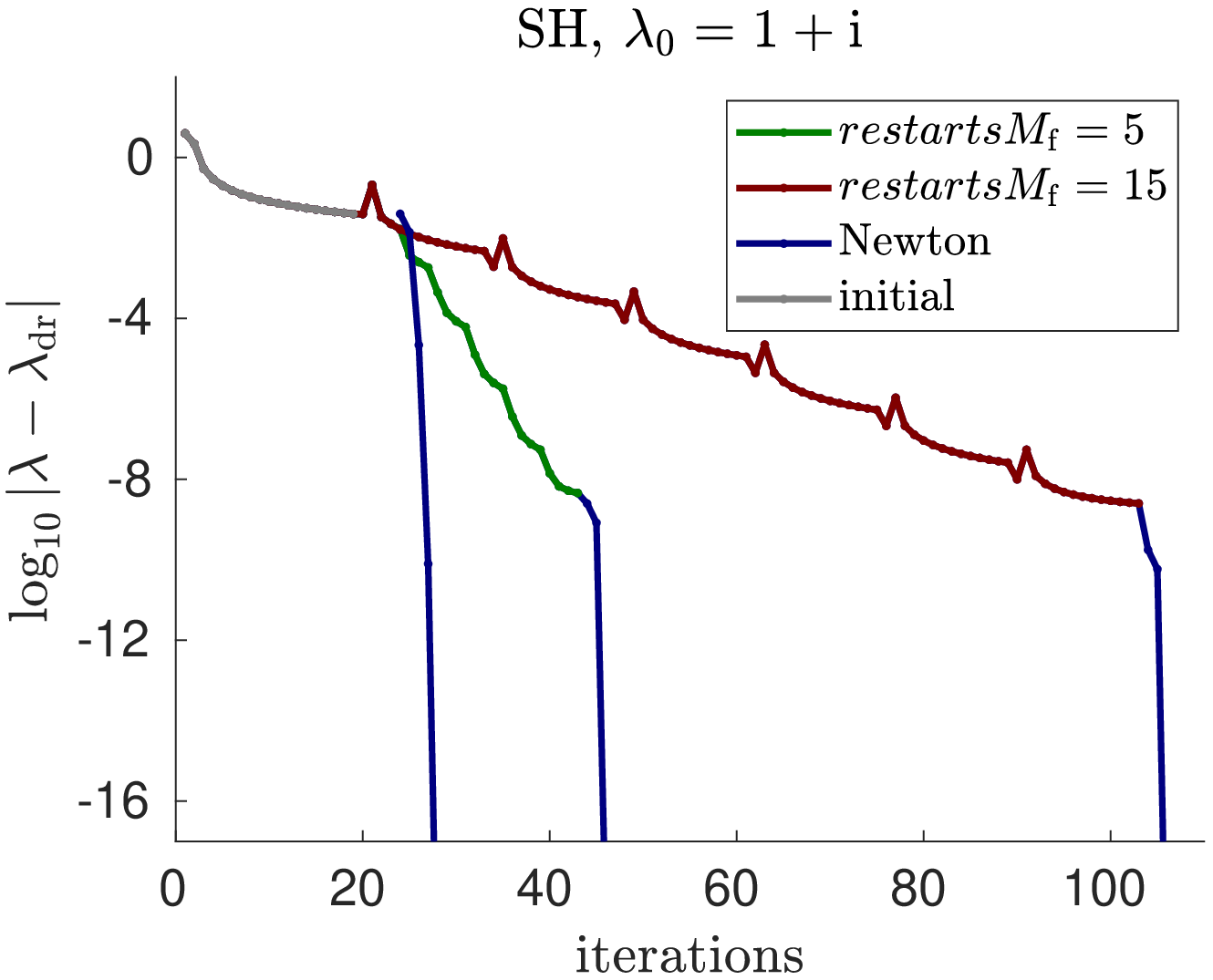}%
\includegraphics[width=0.33\textwidth]{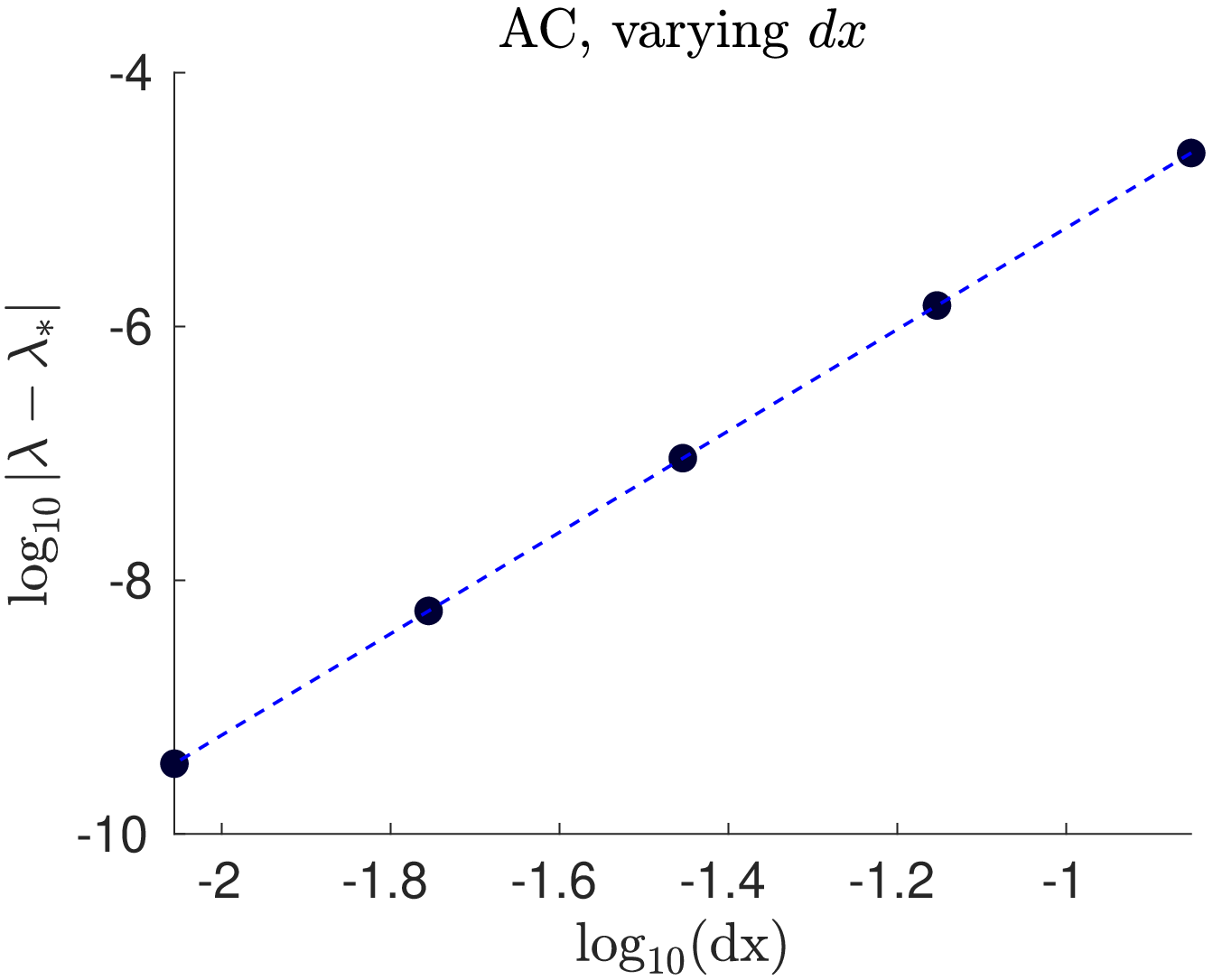}%
\includegraphics[width=0.33\textwidth]{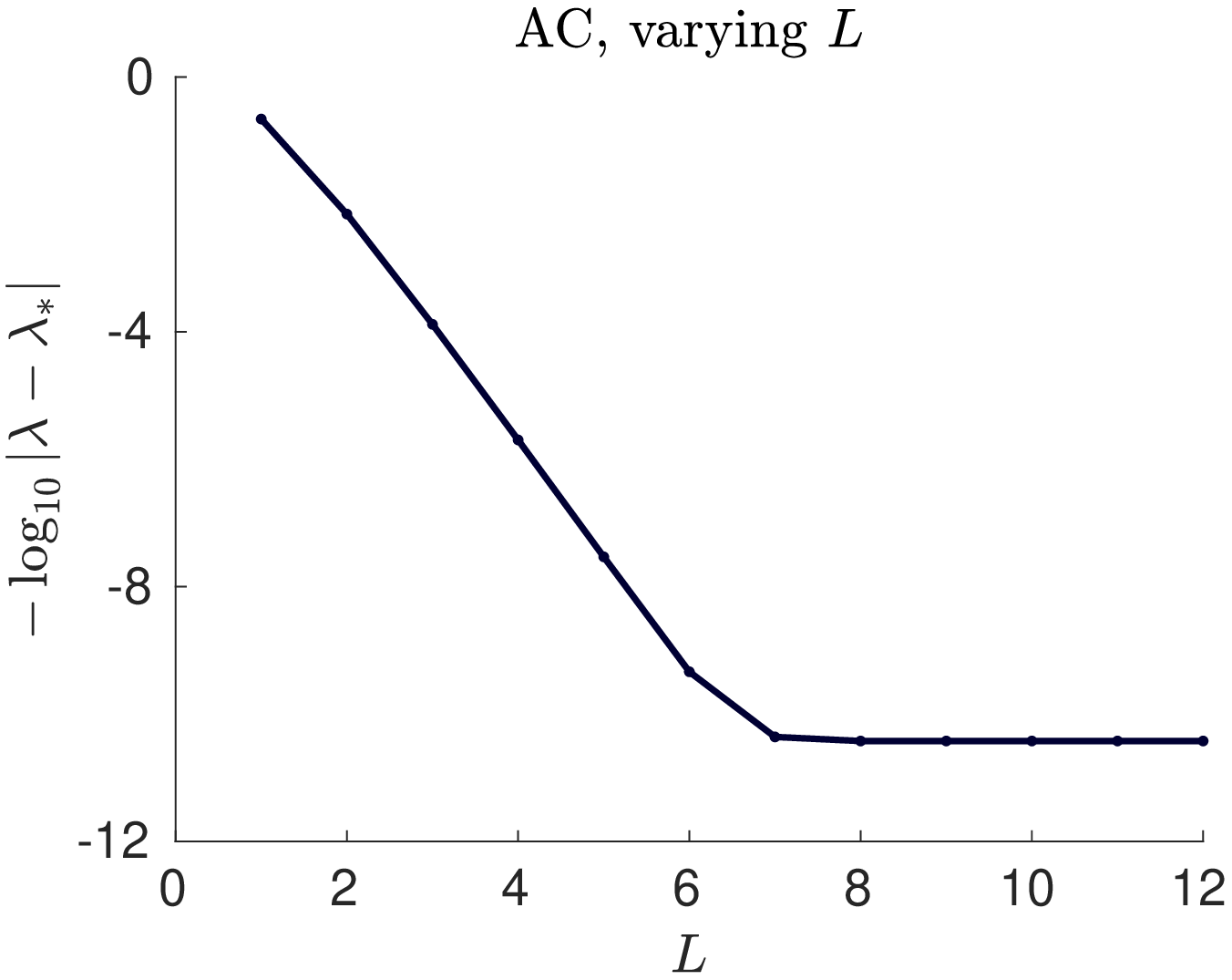}%
\caption{\emph{Left:} Convergence to $\lambda_\mathrm{dr}=0$ in the Swift-Hohenberg (SH) equation \eqref{e:sh} with restarts after 20 initial iterations, $\lambda_0=1+\rmi$; restarts after additional 5 and 15 iterations, respectively, and Newton after just one restart, demonstrating exponential convergence with restarts and practically immediate convergence with Newton for good initial guesses. \emph{Center:} Fourth order convergence in the grid size to the eigenvalue $\lambda_*=0$ with $L=10$ for \eqref{e:ac}. \emph{Right: } Exponential convergence in the domain size $L$ for $dx=0.005$ for  \eqref{e:ac}.  }\label{f:2}
\end{figure}

In our first example, we compute the branch point $\lambda_\mathrm{dr}=0$ associated with the spatial eigenvalue $\nu_\mathrm{dr}=-1$ in 
\begin{equation}\label{e:cd}
w_t=w_{xx}+2w_x+w,
\end{equation}
with unique double root $\lambda_\mathrm{dr}=0$ and associated $\nu_\mathrm{dr}=-1$, and  with starting guess $\lambda_0=1$. Convergence is as expected algebraic with rate $1/k$ but iteration is stable for a very large number iterations, $k\sim 10^4$; see Fig. \ref{f:1}. We find the predicted algebraic convergence with rate $k^{-1}$ from Proposition \ref{p:a3} up to $10^4$ iterates, demonstrating that high-order Taylor expansions can be effective in this context of analytic matrix pencils. Of course, one would in practice restart the computation once sufficient initial accuracy is achieved; see below and Fig. \ref{f:2}. We also confirmed this algebraic rate of convergence in the Swift-Hohenberg equation,
\begin{equation}\label{e:sh}
 w_t=-(\partial_{xx}+1)^2 w,
\end{equation}
with double root  $\lambda_\mathrm{dr}=0$ and associated $\nu_\mathrm{dr}=\rmi$ or $\nu_\mathrm{dr}=-\rmi$, starting value $\lambda_0=1+\rmi$. Convergence is with the predicted rate $k^{-1}$, although $\iota(0)$ has 2-dimensional kernel associated with the two spatial roots $\nu=\pm\rmi$; see Fig. \ref{f:1}, center panel. The Newton method described above indeed identifies both roots. The last example, shown in Fig. \ref{f:1}, right panel, is the linearization at a constant state in the Cahn-Hilliard equation, exhibiting a spinodal decomposition instability. We consider the linearization in a comoving frame such that the double roots $\lambda_\mathrm{dr}=\rmi\omega_\mathrm{dr}$ have zero real part \cite{scheel2017spinodal},
\begin{equation}\label{e:ch}
 w_t=-w_{xxxx}-w_{xx}+c_\mathrm{lin}w_x, \qquad 
 c_\mathrm{lin}=\frac{2}{3\sqrt{6}}\left(2+\sqrt 7 \right)\sqrt{\sqrt 7 -1}, \quad
 \lambda_\mathrm{dr}=\pm\rmi \left(3+\sqrt 7 \right)\sqrt{\frac{2+\sqrt 7}{96}}.
\end{equation}
We also tested convergence for multiple double roots using the Korteweg-De Vries equation 
\begin{equation}\label{e:kdv}
  w_t=w_{xxx}, \qquad \lambda_\mathrm{dr}=0,\ \nu_\mathrm{dr}=0,
\end{equation}
and the beam equation, 
\begin{equation}\label{e:beam}
  w_{tt}=-w_{xxxx}, \qquad \lambda_\mathrm{dr}=0,\ \nu_\mathrm{dr}=0,
\end{equation}
finding the same algebraic convergence rate $k^{-1}$; see Fig. \ref{f:1}, right panel. Convergence to double roots in coupled transport equations from Example \ref{ex:cpw},
\begin{equation}\label{e:ct}
 w^1_t=-w^1_x+\eps w^2,\qquad w^2_t=w^2_x,
\end{equation}
is exponential as expected, since the dispersion relation does not have a branch point at $\lambda_\mathrm{dr}=0$ but rather stable and unstable eigenspaces intersect nontrivially. For $\eps=0$, subspaces do not intersect, the double root disappears. The algorithm picks up this sensitivity through a long transient for small values of $\eps$, before exponential convergence sets in. 

Speed of convergence depends on the distance to the branch point. One therefore would usually first perform a global search for possible instabilities through identifying the closest branch point to an unstable $\lambda_0$. As a second step, one would then try to compute this branch point more precisely through restarting the algorithm with a nearby initial guess as described in \S\ref{s:4} with restarts once increments in the predicted value of $\lambda_\mathrm{dr}$ are small. The result is exponential convergence as demonstrated in Fig. \ref{f:2}, left panel. Typically, one would perform a minimum number of iterations, for instance 5, before repeated restarts since more frequent restarts yield faster convergence. With errors in $\lambda_\mathrm{dr}$ small enough, typically $10^{-3}$, one would switch to a Newton method which will give machine accuracy results within 3 steps.

It is at this point interesting to also return to Example \ref{ex:bp}, $\lambda w=w_{xx}$ on $x>0$ with boundary condition $n_1w+n_2 w_x=0$. Our algorithm identifies (correctly) $\lambda=0$ as a spectral value of $\iota$ regardless of the choice of $n_{1/2}$. Removing this branch point singularity through the choice $\lambda=\gamma^2$ removes the branch singularity and our algorithm finds the spectral values $\gamma=n_1/n_2$, regardless of whether they correspond to eigenvalues, $\gamma>0$, or resonances, $\gamma<0$.

\paragraph{Variable coefficients --- branch points, resonances, and eigenvalues.}

\begin{figure}
\includegraphics[width=0.33\textwidth]{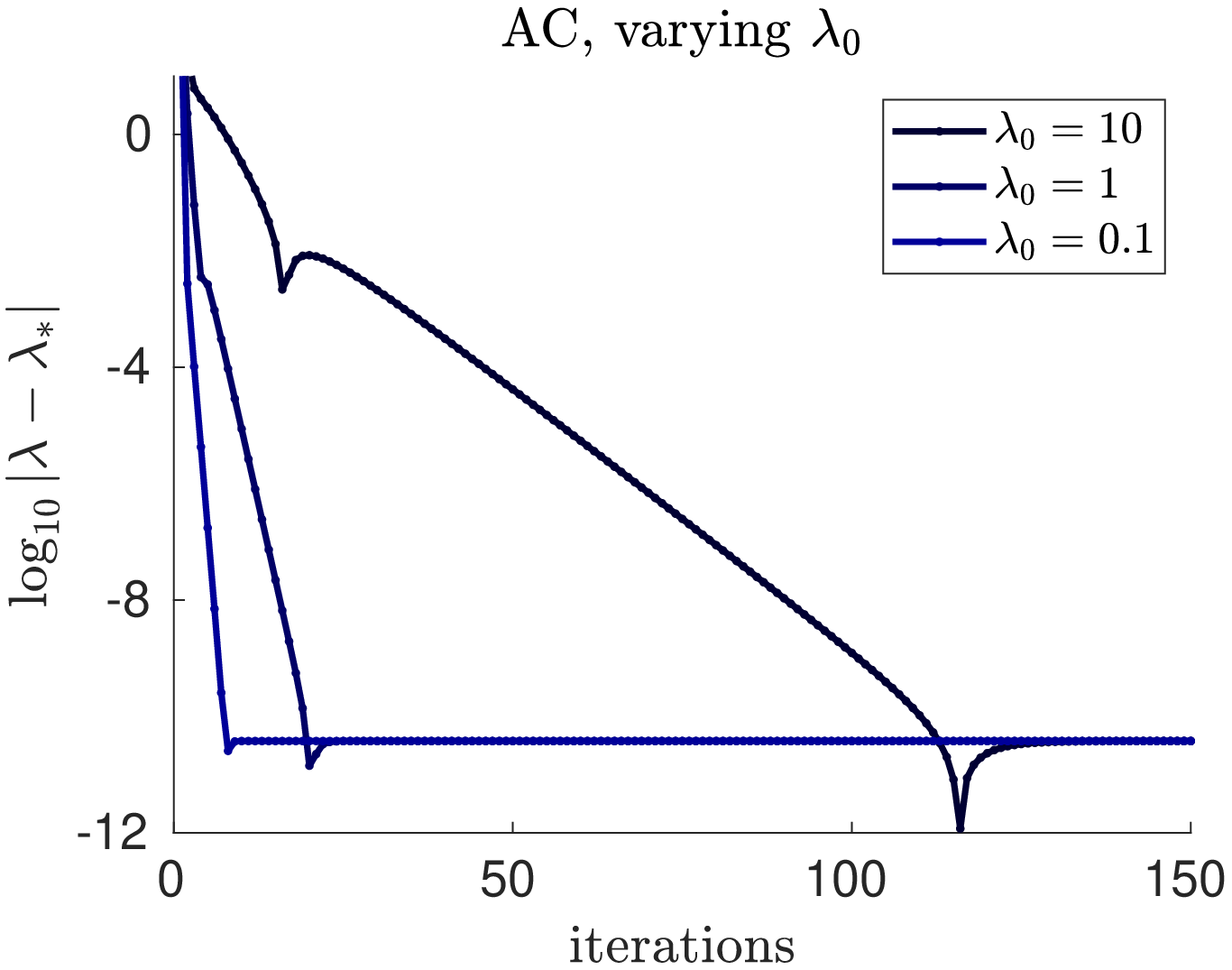}%
\includegraphics[width=0.33\textwidth]{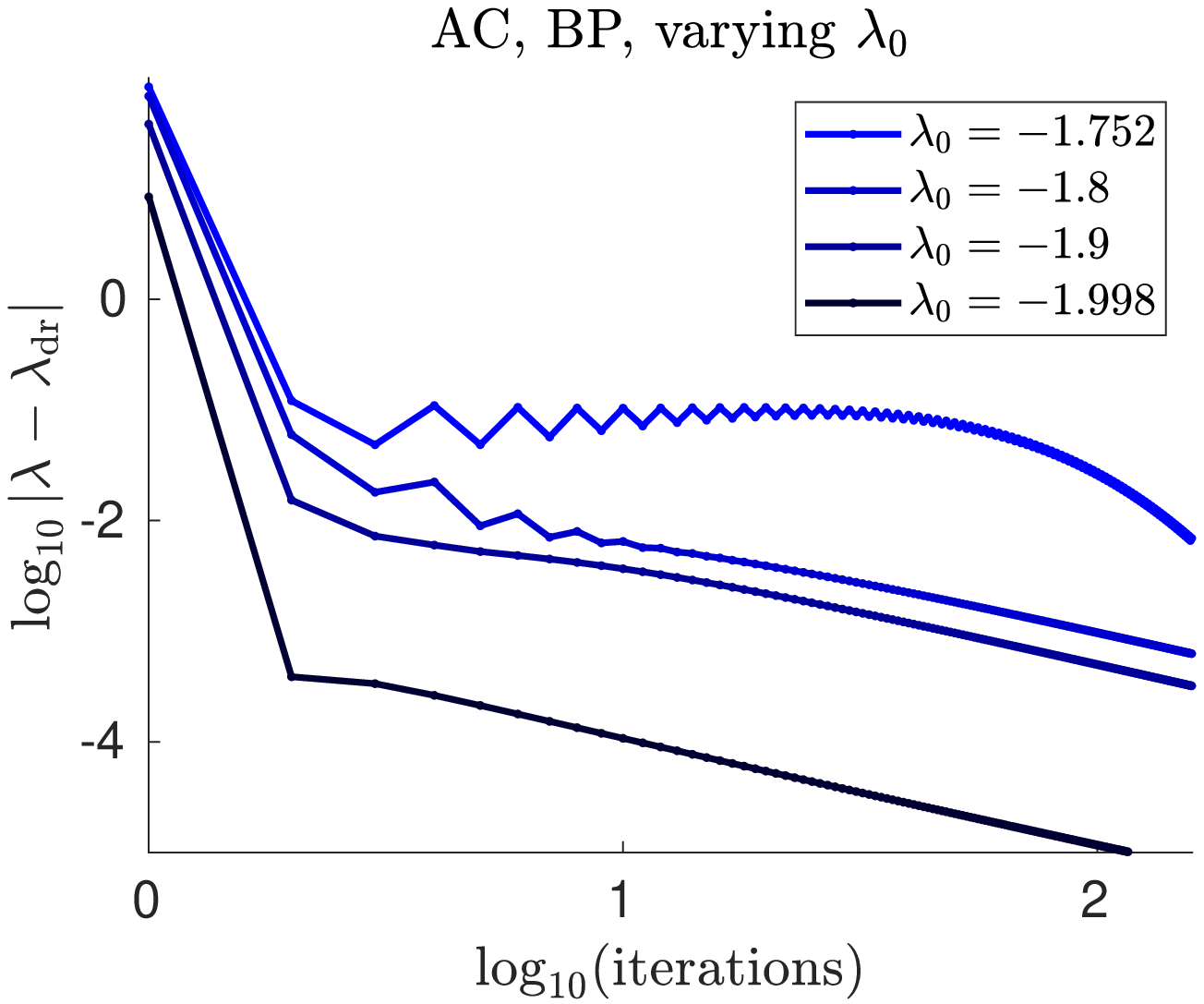}%
\includegraphics[width=0.33\textwidth]{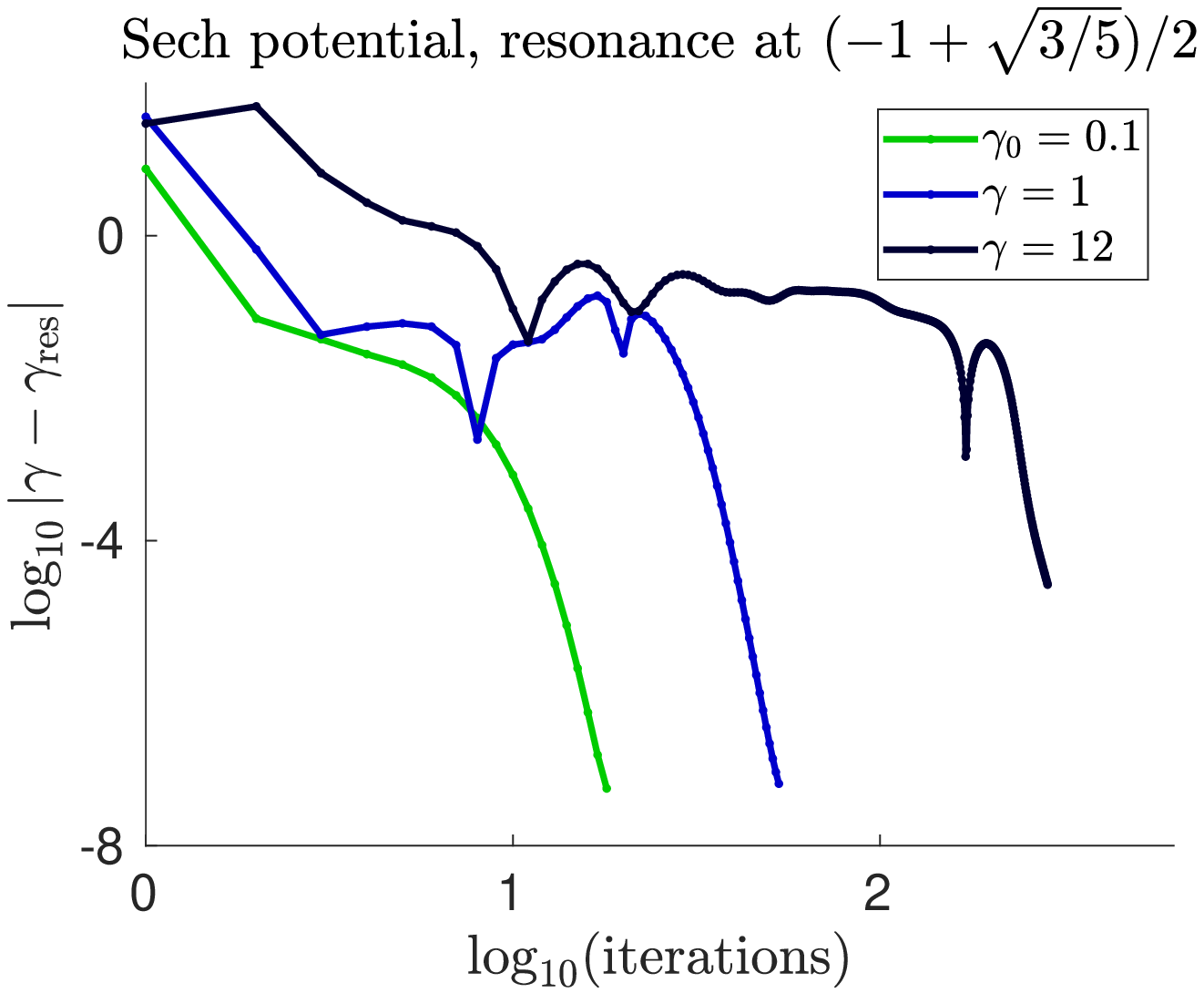}%
\caption{\emph{Left:} Exponential convergence to the eigenvalue $\lambda=0$ in \eqref{e:ac} with convergence rate increasing as $\lambda_0\to 0$. \emph{Center:} Convergence to the branch point $\lambda_\mathrm{dr}=-2$ for different starting values $\lambda_0$; see text for details. \emph{Right:} Convergence to a resonance in \eqref{e:sech} past the domain of analyticity of $\iota$; see text for details.}\label{f:3}
\end{figure}

We illustrate the performance of our algorithm in the case of variable, asymptotically constant coefficients. We start with a 4th-order discretization with grid size $dx$ of the Allen-Cahn layer from Example \ref{ex:1ctd},
\begin{equation}\label{e:ac}
 \lambda w=w_{xx}+(1-3\tanh^2(x/\sqrt{2}))w,
\end{equation}
with eigenvalues at $0$ and $-\frac{3}{2}$, and a branch point at $-2$. The center and right panel in Fig. \ref{f:2} demonstrate 4th order convergence of the compute eigenvalue $\lambda\sim \lambda_*=0$ as $dx$ is decreased in a domain of size $L=10$, and exponential convergence for $dx=0.005$ as $L$ increases. 
Convergence to the eigenvalue is exponential with rate depending on the distance from the eigenvalue (more precisely, the relative distance between the nearest and next-nearest eigenvalue $|\lambda_0-\lambda_1|/|\lambda_0-\lambda_2|)$, with $\lambda_1=0,\lambda_2=-1.75$), which we illustrate in Fig. \ref{f:3}, left panel, with $L=10$ and $dx=0.05$; compare also Proposition \ref{p:pm1} and its proof. Convergence to the branch point $\lambda_\mathrm{dr}=-2$ is algebraic as shown in Fig. \ref{f:3}, center panel; compare also Proposition \ref{p:a3}. However, an initial approach is fast, in particular for starting values close to the branch point, as reflected in \eqref{e:bplamasy}. In fact, restarting the algorithm yields exponential convergence. For starting values close to $-1.75$, branch point and eigenvalue at $\lambda=-1.5$ are at a similar distance and  convergence only sets in after a long transient. We also computed the resonances at $\lambda=-1.5$ and $\lambda=0$ with the same convergence rates, simply exchanging stable and unstable subspaces at $\pm\infty$, confirming the convergence from Proposition \ref{p:pm2}.

Lastly, we present a computation of resonances in 
\begin{equation}\label{e:sech}
 \lambda w = w_{xx}+F_0 \mathrm{sech}^2(x) w,\qquad F_0=-0.1,\quad \gamma_\mathrm{res}=-\frac{1}{2} \sqrt{F_0+\frac{1}{4}},\quad \lambda=\gamma^2.
\end{equation}
Writing $\lambda=\gamma^2$, removes the branch point at $\lambda=0$ and allows for detection of the resonance closest to $\gamma_\mathrm{res}$. We use $F_0=-1/10$, which gives $\gamma_\mathrm{res}=(-1+\sqrt{3/5})/2  \sim -0.1127$. The stable subspace at $\gamma_0>0$ is given by $(1,\gamma_0)^T$. Writing eigenspaces as graphs over this subspace yields a pole at $\gamma=-1/\gamma_0$. In particular, for $\gamma_0=12$, the series expansion of the boundary condition has a pole at $\gamma=-1/12\sim -0.0833$, between $\gamma_0$ and $\gamma_\mathrm{res}$, so that $\gamma_\mathrm{res}$ is not located within the radius of convergence of $\iota$ when choosing this initial value. Fig. \ref{f:3}, right panel, demonstrates convergence in this situation as predicted by Proposition \ref{p:pm2}. Convergence is slow and can again be accelerated using restarts, as is clear from the rates of convergence for initial guesses closer to $\gamma_\mathrm{res}$.

Computation times are all less than 10 seconds, with the exception of the example in Fig. \ref{f:1}, left panel, where a very large number of iterations was performed and a very high order of the Taylor expansion needs to be precomputed, leading to computation times of roughly 3 minutes on a laptop.

% \includegraphics[width=0.33\textwidth]{}%
% \includegraphics[width=0.33\textwidth]{}%
% \includegraphics[width=0.33\textwidth]{}%
% \includegraphics[width=0.33\textwidth]{}%

% 
% ac_dx.eps        conv_diff_no_restart.eps  sh_ch_no_restart.eps
% ac_lam_0_bp.eps  copy_figures              sh_w_restart.eps
% ac_lam_0.eps     kdv_wave_no_restart.eps
% ac_L.eps         sech_res.eps

\paragraph{Large problems.}

Elliptic problems in cylindrical domains $(x,y)\in \R\times\Omega$ yield, after discretization in the $y$-direction, high-dimensional problems of the form \eqref{e:twlin}, $N\gg 1$. We demonstrate that the methods here are capable of treating such problems with the example of a Schr\"odinger eigenvalue problem with a localized potential trap and absorbing boundary conditions,
\begin{equation}\label{e:sch_str}
 w_{xx}+w_{yy}+\varepsilon V(x,y)w=\lambda w, \qquad (x,y)\in \R\times (-\pi,\pi),\qquad  w(x,\pm\pi)=0, \qquad V(x,y)=\mathrm{sech}^2\,\left(\sqrt{x^2+y^2}\right).
\end{equation}
One readily finds the essential spectrum at $(-\infty,-\frac{1}{4}]$ terminating in a branch point at $-\frac{1}{4}$, which we resolve by considering the problem on the Riemann surface with new eigenvalue parameter $\gamma=\sqrt{\lambda+\frac{1}{4}}$ and branch cut of the square root at the negative real line, which gives
\begin{equation}\label{e:sch_srt2}
 w_{xx}+w_{yy}+\varepsilon V(x,y)w=\left(\gamma^2-\frac{1}{4}\right) w, \qquad (x,y)\in \R\times (-\pi,\pi),\qquad  w(x,\pm\pi)=0.
\end{equation}
For $\varepsilon\gtrsim 0$, the eigenvalue problem possesses a unique eigenvalue $\lambda_*(\varepsilon)$, with expansion
\begin{equation}\label{e:evex}
 \lambda_*(\varepsilon)=\gamma_*(\varepsilon)^2-\frac{1}{4},\qquad \gamma_*(\varepsilon)=\gamma_1 \varepsilon+\rmO(\varepsilon^2),\qquad \gamma_1=\frac{1}{2\pi}\int_{x,y}V(x,y)\cos^2(y/2)=  0.567402\ldots,
\end{equation}
using \cite{SIMON1976279} with technical adaptations as in \cite{MR3283552}. Eigenfunctions have asymptotics $u(x,y)\sim \cos(y/2)\rme^{-\nu |x|}$ for $|x|\to\infty$, with $\nu=\gamma=\sqrt{\lambda+\frac{1}{4}}$.

For $\varepsilon\lesssim 0$, the eigenvalue changes into a resonance pole at $\lambda_*(\varepsilon)$, with the same expansion \eqref{e:evex}. The eigenfunction exhibits asymptotic growth $u(x,y)\sim \cos(y/2)\rme^{\nu |x|}+\rmO(\rme^{-\eta|x|})$ for $|x|\to\infty$, with $\nu=\gamma=\sqrt{\lambda+\frac{1}{4}}$, and $\eta=\rmO(1)$ as $\eps\to 0$.

Truncating the unbounded strip to $(x,y)\in (-L,L)\times(-\pi,\pi)$ with say Dirichlet boundary conditions at $x=\pm L$ yields truncation errors for eigenvalues, $\varepsilon>0$, of order $\rme^{-2\nu L}$, thus requiring $L\gg 1/\varepsilon$. The essential spectrum breaks into clusters of eigenvalues with gaps $\rmO(L^{-2})$ starting at $\lambda=-1/4$; see \cite{ssabs}. Resonances $\lambda_*(\varepsilon)$,  $\varepsilon<0$, cannot be found easily in such truncations.

In our approach, we rewrite \eqref{e:sch_srt2} in the form 
\begin{equation}\label{e:sch_str_sys}
 u_x= A(x;\lambda) u,\qquad  A(x;\lambda)=
 \begin{pmatrix}
 0&1\\
 -\partial_{yy}-\varepsilon V(x,y)+\gamma^2-\frac{1}{4} & 0
 \end{pmatrix},
\end{equation}
and discretize $\partial_{yy}$, using fourth order centered finite differences and Dirichlet boundary conditions. We tested spatial discretizations at $\eps=1.5$, finding that the $y$-discretization error is well below $10^{-8}$ with $N_y=300$, based on a reference eigenvalue with $N_y=450$. We used comparable discretization in $x$, so that $dx\sim dy \sim 0.02$. 
\begin{figure}
\includegraphics[width=0.33\textwidth]{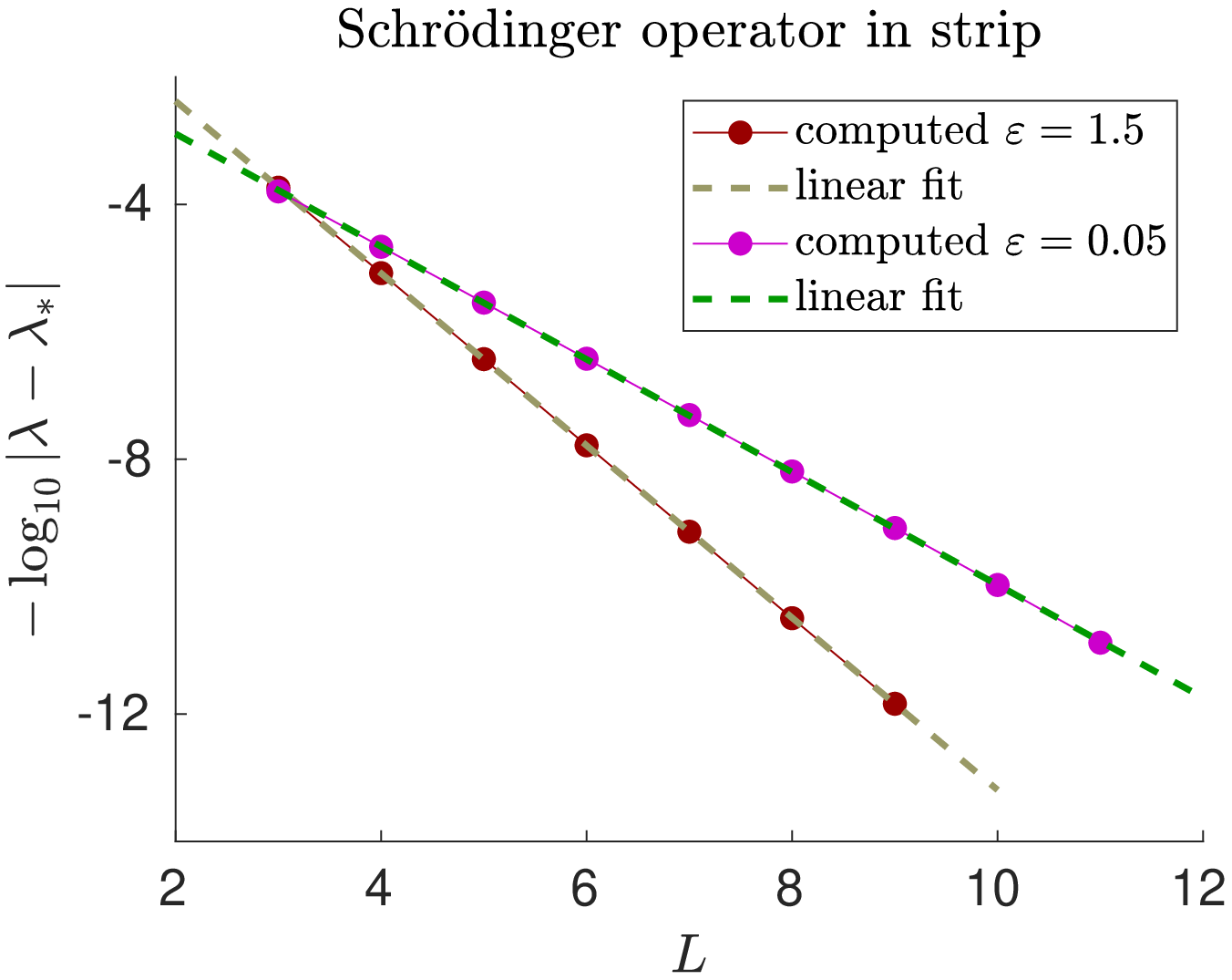}%
\includegraphics[width=0.33\textwidth]{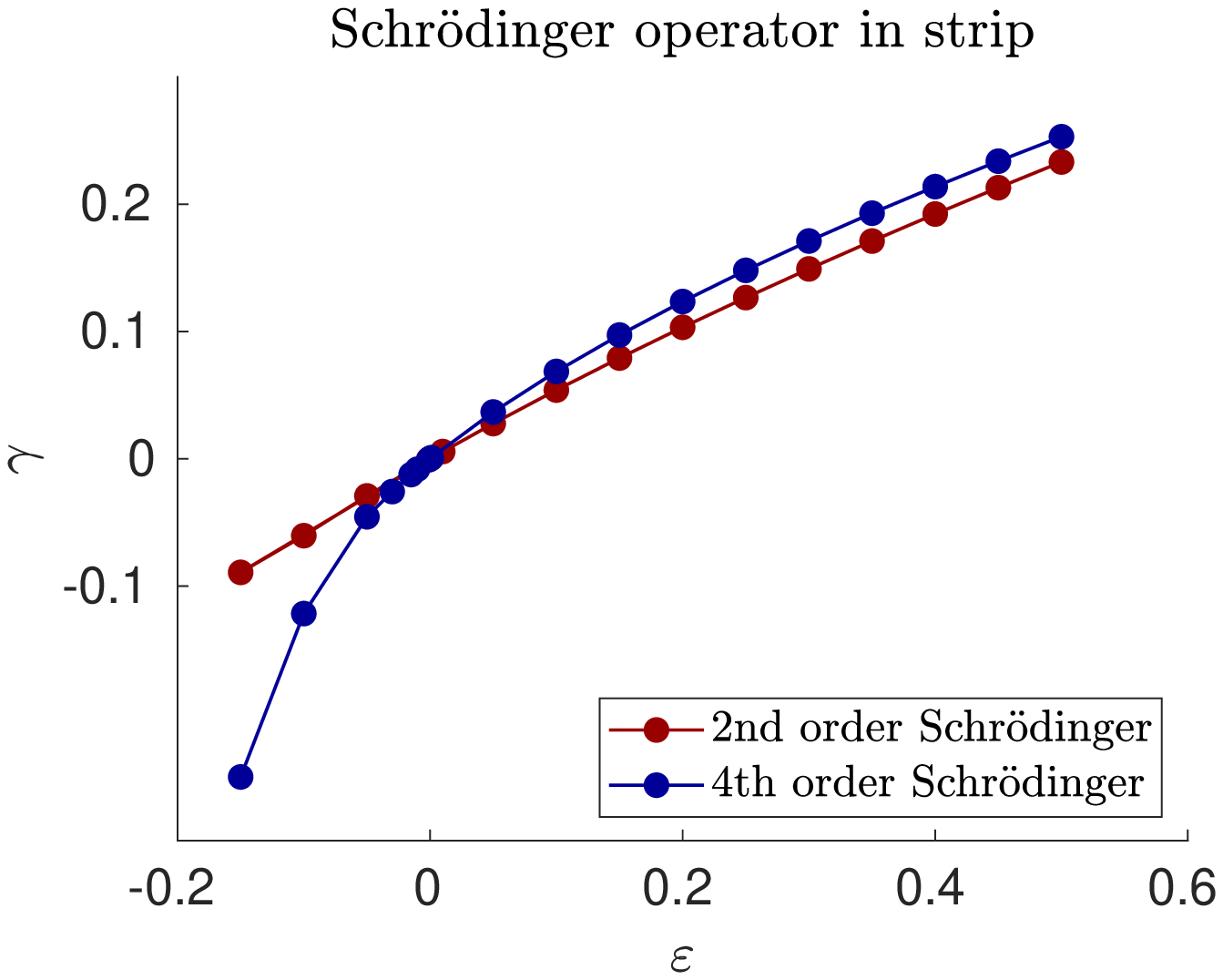}%
\includegraphics[width=0.33\textwidth]{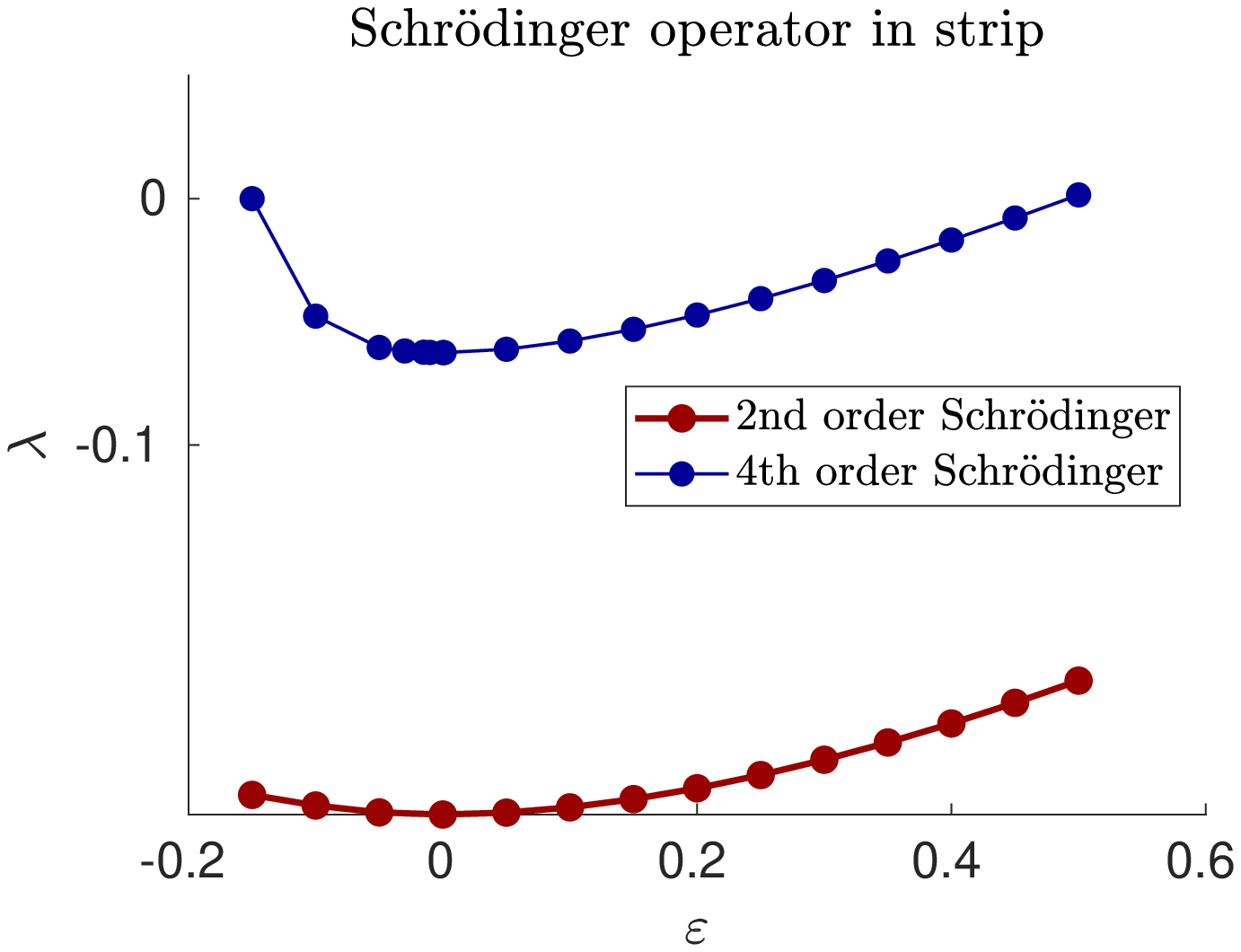}%
 \caption{\emph{Left:} Exponential convergence of computed eigenvalues in \eqref{e:sch_str} with domain size $L$ to limits $\lambda_*=0.076763657389$ ($\eps=1.5$) and $\lambda_*=-0.24923674$ ($\eps=0.05$); slopes from linear fit $-3.11$ ($\eps=1.5$) and $-2.04$ ($\eps=0.05$) correspond well to theoretical predictions $-3.14$ and $-2.054$, respectively. \emph{Center and left:} Smooth continuation of eigenvalues $\lambda=\gamma^2-1/4$ into resonance poles as functions of $\varepsilon$, with starting value $\gamma=0.1$. Best linear approximation at $\eps=0$ from $\eps=10^{-3}$ gives  $\gamma\sim 0.56709\cdot\eps$ which compares well with \eqref{e:evex}. Also shown are results for a fourth order Schr\"odinger equation with the same potential, \eqref{e:4th}, with linear expansion  $\gamma\sim 0.802428\cdot\eps$, $\lambda=\gamma^2-\frac{1}{16}$, again in excellent agreement with theory \eqref{e:4thexp}.  }\label{f:strip}
\end{figure}
%     -8.024284733314771e-04
%
Determinants of $\iota^0$ evaluated to infinity for even moderate grids and any attempt at finding eigenvalues using winding number arguments for determinants would likely require renormalizations, using for instance  Fredholm determinants, for which however numerical computations are not well developed \cite{MR2600548}.

The asymptotic boundary conditions used in our formulation imply convergence with a uniform rate when $\eps\sim 0$ as $L\to\infty$. Boundary conditions are accurate with rate $\rme^{-2L}$ given by the convergence rate of $V$, leading to a predicted error  $\rme^{(-2-2\gamma)L}$, matching well the numerically observed error shown in Fig. \ref{f:strip} (left panel). 

We used $L=8$, $N_y=300$, $N_x=800$, $dy\sim dx=0.02$, to compute the eigenvalue for different values of $\eps$. For starting values $\lambda_0=0.1$ we found convergence with errors of order $10^{-6}$ within 40 primary iterations. Clearly, a continuation approach would be more effective for computing the resulting curves of eigenvalues and resonances shown in Fig. \ref{f:strip} (center and left panel) --- the rapid convergence for fixed starting values however  demonstrates the broader applicability of our approach. Asymptotics near $\eps=0$ agree very well with the prediction from \eqref{e:evex}; see Fig. \ref{f:strip} for comparison.

In the slightly harder problem with fourth order dispersion,
\begin{equation}\label{e:4th}
-\left( \partial_{xx}+\partial_{yy}\right)^2+\varepsilon V(x,y)w=\lambda w, \qquad (x,y)\in \R\times (-\pi,\pi),\qquad  w(x,\pm\pi)=0, \qquad V(x,y)=\mathrm{sech}^2\,\left(\sqrt{x^2+y^2}\right).
\end{equation}
one can mimic the theoretical prediction near the edge of the essential spectrum $\lambda_\mathrm{bp}=-\frac{1}{16}$ and find the expansion for eigenvalue and resonance at $\eps=0$,
\begin{equation}\label{e:4thexp}
\gamma=\sqrt{2}\gamma_1\eps+\rmO(\eps^2),\qquad \sqrt{2}\gamma_1=0.802428\ldots ,
\end{equation}
with $\gamma_1$ from \eqref{e:evex}.

We also computed eigenvalues for the potential $V(x,y)=\frac{1}{2}\,\mathrm{sech}^2(\frac{x}{2})$, in particular $\lambda=0$ with eigenfunction $\cos(\frac{y}{2})\,\mathrm{sech}\,(\frac{x}{2})$, confirming the convergence rates documented above when an explicit eigenvalue is known. 

Computation times are several minutes on a laptop. Memory requirements limit the system size since matrices are full in the index for the $y$-component. A sparse approximation of $\iota$, even for constant-in-$x$ problems would be desirable.

\paragraph{Spreading speeds.} Localized disturbances of an unstable state grow temporally and spread spatially. The spatial spreading can be captured via the study of pointwise instabilities in comoving frames; see \cite{holzerscheel14} for background. Using the algorithms above, one would compute branch points in a constant-coefficient problem 
\[
 \lambda w=\mathcal{P}(\partial_x)w, \qquad \text{ or } \quad u_x=A(\lambda)u.
\]
One would then track double roots $\lambda_\mathrm{dr}$ with associated spatial exponent $\nu_\mathrm{dr}$ using numerical continuation as a function of $c$ in 
\begin{equation}\label{e:comov}
 \lambda w=\mathcal{P}(\partial_x)w+cw_x, \qquad \text{ or }\quad u_x=\tilde{A}(\lambda,c)u .
\end{equation}
Increasing $c$, one tracks $\lambda_\mathrm{dr}(c)$ and finds the largest value $c_\mathrm{lin}$ of $c$ so that $\Re\lambda_\mathrm{dr}(c)=0$. One would then, for this specific value of $c$ verify that there are no unstable double roots, leaving open however the possibility of instabilities for yet larger values of $c$. 

We mention here a more direct method that yields directly critical values $c_\mathrm{lin}$ in the case where the associated branch point $\lambda_\mathrm{dr}$ is real.  One therefore simply considers \eqref{e:comov} with $\lambda=0$,
\begin{equation}
 u_x=\tilde{A}(0,c)u,
\end{equation}
as a nonlinear eigenvalue problem in $c$! ``Eigenvalues'' $c$ correspond to values of $c$ where pointwise growth is neutral, $\lambda_\mathrm{dr}=0$, and thus yield all candidates for linear spreading speeds, with the largest one typically being most relevant. We verified numerically that this algorithm performs very well in the extended Fisher-KPP equation,
\[
 w_t=-\eps^2 w_{xxxx}+w_{xx}+w-w^3, \text{ with linearization } w_t=-\eps w_{xxxx}+w_{xx}+w, 
\]
and spreading speeds 
\[
 c_\mathrm{lin}=\frac{1}{9} \sqrt{\frac{6-6 \sqrt{1-12
   {\eps}^2}}{{\eps}^2}} \left(\sqrt{1-12
   {\eps}^2}+2\right), \qquad \text{for } \eps^2<\frac{1}{12}.
\]
We note that spreading speeds may be (and indeed are in this example for $\eps^2>\frac{1}{12}$)  associated with complex values $\lambda_\mathrm{dr}=\rmi\omega_\mathrm{dr}$, which are not detected by this procedure. The algorithm rather yields complex speeds $c_\mathrm{lin}$ which do not appear to be relevant to the stability problem.

\section{Summary and outlook}\label{s:6}
We proposed an inverse power method as a versatile tool to locate spectral values of differential operators on the real line. The method identifies all singularities of the pointwise Green's function, including eigenvalues, resonances, and branch points, finding in particular the closest singularity to a given reference point $\lambda_*$. Pointwise methods have been used mostly in connection with the Evans function, effectively taking determinants. We hope that our view point provides a robust alternative to such determinant-based methods and will then prove useful particularly in large systems. 
In future work, we plan to investigate effective strategies for large systems, when bases for stable and unstable subspaces yield full matrices for $\iota$, and the case of periodic coefficients. On the other hand, it appears to be difficult to adapt this formalism to yield spreading speeds also in the oscillatory case $\omega_\mathrm{dr}\neq 0$, and to multi-dimensional problems. Similarly, the pointwise formulation adapted here relies strongly on a ``local'' formulation in $x$, excluding to some extend spatially nonlocal coupling that does not permit a formulation as a first-order spatial ODE through linearization of the matrix pencil in $\partial_x$; see however \cite{MR3283552,MR3803149,MR4309433} for techniques that recover ``pointwise'' descriptions in this nonlocal setting. Similarly, effective computational tools to analyze multi-dimensional problems in this pointwise context do not appear to be available; see for instance \cite{MR1140700} for a discussion of pointwise instabilities in constant-coefficient, multi-dimensional problems.

\paragraph{Acknowledgments.} The author acknowledges partial support through grants NSF DMS-1907391 and DMS-2205663.
\paragraph{Code.} Code used for the computations in the examples is available at the repository \textsc{https://github.com/arnd-scheel/nonlinear-eigenvalue}

\def\cprime{$'$}

% 
% \bibliographystyle{abbrv}
% \bibliography{nl}

\begin{thebibliography}{10}

\bibitem{MR1068805}
J.~Alexander, R.~Gardner, and C.~Jones.
\newblock A topological invariant arising in the stability analysis of
  travelling waves.
\newblock {\em J. Reine Angew. Math.}, 410:167--212, 1990.

\bibitem{https://doi.org/10.48550/arxiv.2012.06443}
M.~Avery and A.~Scheel.
\newblock Universal selection of pulled fronts.
\newblock {\em Communications of the AMS, to appear}, 2022.

\bibitem{stablab}
B.~Barker, J.~Humpherys, J.~Lytle, and K.~Zumbrun.
\newblock Stablab: A matlab-based numerical library for evans function
  computation.
\newblock {\em Available in the github repository, nonlinear-waves/stablab.},
  2015.

\bibitem{MR3789546}
M.~Beck, G.~Cox, C.~Jones, Y.~Latushkin, K.~McQuighan, and A.~Sukhtayev.
\newblock Instability of pulses in gradient reaction-diffusion systems: a
  symplectic approach.
\newblock {\em Philos. Trans. Roy. Soc. A}, 376(2117):20170187, 20, 2018.

\bibitem{MR3157977}
W.-J. Beyn, Y.~Latushkin, and J.~Rottmann-Matthes.
\newblock Finding eigenvalues of holomorphic {F}redholm operator pencils using
  boundary value problems and contour integrals.
\newblock {\em Integral Equations Operator Theory}, 78(2):155--211, 2014.

\bibitem{MR3144797}
D.~Bindel and A.~Hood.
\newblock Localization theorems for nonlinear eigenvalue problems.
\newblock {\em SIAM J. Matrix Anal. Appl.}, 34(4):1728--1749, 2013.

\bibitem{MR2600548}
F.~Bornemann.
\newblock On the numerical evaluation of {F}redholm determinants.
\newblock {\em Math. Comp.}, 79(270):871--915, 2010.

\bibitem{MR1140700}
L.~Brevdo.
\newblock Three-dimensional absolute and convective instabilities, and
  spatially amplifying waves in parallel shear flows.
\newblock {\em Z. Angew. Math. Phys.}, 42(6):911--942, 1991.

\bibitem{brevdo_linear_1999}
L.~Brevdo, P.~Laure, F.~Dias, and T.~J. Bridges.
\newblock Linear pulse structure and signalling in a film flow on an inclined
  plane.
\newblock {\em Journal of Fluid Mechanics}, 396:37--71, 1999.

\bibitem{doi:10.1146/annurev.fluid.37.061903.175810}
J.-M. Chomaz.
\newblock Global instabilities in spatially developing flows: Non-normality and
  nonlinearity.
\newblock {\em Annual Review of Fluid Mechanics}, 37(1):357--392, 2005.

\bibitem{PhysRevE.105.014602}
C.~del Junco, A.~Estevez-Torres, and A.~Maitra.
\newblock Front speed and pattern selection of a propagating chemical front in
  an active fluid.
\newblock {\em Phys. Rev. E}, 105:014602, Jan 2022.

\bibitem{MR1878337}
A.~Doelman, R.~A. Gardner, and T.~J. Kaper.
\newblock A stability index analysis of 1-{D} patterns of the {G}ray-{S}cott
  model.
\newblock {\em Mem. Amer. Math. Soc.}, 155(737):xii+64, 2002.

\bibitem{faye17}
G.~Faye, M.~Holzer, and A.~Scheel.
\newblock Linear spreading speeds from nonlinear resonant interaction.
\newblock {\em Nonlinearity}, 30(6):2403--2442, may 2017.

\bibitem{MR3283552}
G.~Faye and A.~Scheel.
\newblock Fredholm properties of nonlocal differential operators via spectral
  flow.
\newblock {\em Indiana Univ. Math. J.}, 63(5):1311--1348, 2014.

\bibitem{MR3803149}
G.~Faye and A.~Scheel.
\newblock Center manifolds without a phase space.
\newblock {\em Trans. Amer. Math. Soc.}, 370(8):5843--5885, 2018.

\bibitem{fiedler03}
B.~Fiedler and A.~Scheel.
\newblock Spatio-temporal dynamics of reaction-diffusion patterns.
\newblock In {\em Trends in nonlinear analysis}, pages 23--152. Springer,
  Berlin, 2003.

\bibitem{PhysRevB.83.064113}
P.~K. Galenko and K.~R. Elder.
\newblock Marginal stability analysis of the phase field crystal model in one
  spatial dimension.
\newblock {\em Phys. Rev. B}, 83:064113, Feb 2011.

\bibitem{gardner98}
R.~A. Gardner and K.~Zumbrun.
\newblock The gap lemma and geometric criteria for instability of viscous shock
  profiles.
\newblock {\em Comm. Pure Appl. Math.}, 51(7):797--855, 1998.

\bibitem{MR2350362}
F.~Gesztesy, Y.~Latushkin, and K.~A. Makarov.
\newblock Evans functions, {J}ost functions, and {F}redholm determinants.
\newblock {\em Arch. Ration. Mech. Anal.}, 186(3):361--421, 2007.

\bibitem{gohberg}
I.~Gohberg, P.~Lancaster, and L.~Rodman.
\newblock {\em Invariant subspaces of matrices with applications}, volume~51 of
  {\em Classics in Applied Mathematics}.
\newblock Society for Industrial and Applied Mathematics (SIAM), Philadelphia,
  PA, 2006.
\newblock Reprint of the 1986 original.

\bibitem{gutteltisseur}
S.~G\"{u}ttel and F.~Tisseur.
\newblock The nonlinear eigenvalue problem.
\newblock {\em Acta Numer.}, 26:1--94, 2017.

\bibitem{holzerscheel14}
M.~Holzer and A.~Scheel.
\newblock Criteria for pointwise growth and their role in invasion processes.
\newblock {\em J. Nonlinear Sci.}, 24(4):661--709, 2014.

\bibitem{MR3413592}
J.~Humpherys and J.~Lytle.
\newblock Root following in {E}vans function computation.
\newblock {\em SIAM J. Numer. Anal.}, 53(5):2329--2346, 2015.

\bibitem{MR2221065}
J.~Humpherys, B.~Sandstede, and K.~Zumbrun.
\newblock Efficient computation of analytic bases in {E}vans function analysis
  of large systems.
\newblock {\em Numer. Math.}, 103(4):631--642, 2006.

\bibitem{MR2253406}
J.~Humpherys and K.~Zumbrun.
\newblock An efficient shooting algorithm for {E}vans function calculations in
  large systems.
\newblock {\em Phys. D}, 220(2):116--126, 2006.

\bibitem{MR3100266}
T.~Kapitula and K.~Promislow.
\newblock {\em Spectral and dynamical stability of nonlinear waves}, volume 185
  of {\em Applied Mathematical Sciences}.
\newblock Springer, New York, 2013.
\newblock With a foreword by Christopher K. R. T. Jones.

\bibitem{kapitula98}
T.~Kapitula and B.~Sandstede.
\newblock Stability of bright solitary-wave solutions to perturbed nonlinear
  {S}chr\"odinger equations.
\newblock {\em Phys. D}, 124(1-3):58--103, 1998.

\bibitem{MR1897705}
T.~Kapitula and B.~Sandstede.
\newblock Edge bifurcations for near integrable systems via {E}vans function
  techniques.
\newblock {\em SIAM J. Math. Anal.}, 33(5):1117--1143, 2002.

\bibitem{kato}
T.~Kato.
\newblock {\em Perturbation theory for linear operators}.
\newblock Grundlehren der Mathematischen Wissenschaften, Band 132.
  Springer-Verlag, Berlin-New York, second edition, 1976.

\bibitem{MR1759902}
G.~J. Lord, D.~Peterhof, B.~Sandstede, and A.~Scheel.
\newblock Numerical computation of solitary waves in infinite cylindrical
  domains.
\newblock {\em SIAM J. Numer. Anal.}, 37(5):1420--1454, 2000.

\bibitem{mennicken}
R.~Mennicken and M.~M\"{o}ller.
\newblock {\em Non-self-adjoint boundary eigenvalue problems}, volume 192 of
  {\em North-Holland Mathematics Studies}.
\newblock North-Holland Publishing Co., Amsterdam, 2003.

\bibitem{palmer}
K.~J. Palmer.
\newblock Exponential dichotomies and {F}redholm operators.
\newblock {\em Proc. Amer. Math. Soc.}, 104(1):149--156, 1988.

\bibitem{MR1177566}
R.~L. Pego and M.~I. Weinstein.
\newblock Eigenvalues, and instabilities of solitary waves.
\newblock {\em Philos. Trans. Roy. Soc. London Ser. A}, 340(1656):47--94, 1992.

\bibitem{rademacher07}
J.~D.~M. Rademacher, B.~Sandstede, and A.~Scheel.
\newblock Computing absolute and essential spectra using continuation.
\newblock {\em Phys. D}, 229(2):166--183, 2007.

\bibitem{sandstede02}
B.~Sandstede.
\newblock Stability of travelling waves.
\newblock In {\em Handbook of dynamical systems, {V}ol. 2}, pages 983--1055.
  North-Holland, Amsterdam, 2002.

\bibitem{ssabs}
B.~Sandstede and A.~Scheel.
\newblock Absolute and convective instabilities of waves on unbounded and large
  bounded domains.
\newblock {\em Phys. D}, 145(3-4):233--277, 2000.

\bibitem{ssmodstab}
B.~Sandstede and A.~Scheel.
\newblock On the structure of spectra of modulated travelling waves.
\newblock {\em Math. Nachr.}, 232:39--93, 2001.

\bibitem{sandstede04}
B.~Sandstede and A.~Scheel.
\newblock Evans function and blow-up methods in critical eigenvalue problems.
\newblock {\em Discrete Contin. Dyn. Syst.}, 10(4):941--964, 2004.

\bibitem{ssmorse}
B.~Sandstede and A.~Scheel.
\newblock Relative {M}orse indices, {F}redholm indices, and group velocities.
\newblock {\em Discrete Contin. Dyn. Syst.}, 20(1):139--158, 2008.

\bibitem{scheel2017spinodal}
A.~Scheel.
\newblock Spinodal decomposition and coarsening fronts in the cahn--hilliard
  equation.
\newblock {\em Journal of Dynamics and Differential Equations}, 29(2):431--464,
  2017.

\bibitem{MR4309433}
W.~M. Schouten-Straatman and H.~J. Hupkes.
\newblock Exponential dichotomies for nonlocal differential operators with
  infinite range interactions.
\newblock {\em J. Differential Equations}, 301:353--427, 2021.

\bibitem{doi:10.1073/pnas.1420171112}
J.~A. Sherratt.
\newblock Using wavelength and slope to infer the historical origin of semiarid
  vegetation bands.
\newblock {\em Proceedings of the National Academy of Sciences},
  112(14):4202--4207, 2015.

\bibitem{shubin}
M.~A. Shubin.
\newblock On holomorphic families of subspaces of a {B}anach space.
\newblock {\em Integral Equations Operator Theory}, 2(3):407--420, 1979.

\bibitem{SIMON1976279}
B.~Simon.
\newblock The bound state of weakly coupled schrödinger operators in one and
  two dimensions.
\newblock {\em Annals of Physics}, 97(2):279--288, 1976.

\bibitem{MR2183609}
S.~A. Suslov.
\newblock Numerical aspects of searching convective/absolute instability
  transition.
\newblock {\em J. Comput. Phys.}, 212(1):188--217, 2006.

\bibitem{trofimov}
V.~P. Trofimov.
\newblock The root subspaces of operators that depend analytically on a
  parameter.
\newblock {\em Mat. Issled.}, 3(vyp. 3 (9)):117--125, 1968.

\bibitem{vansaarloos03}
W.~van Saarloos.
\newblock Front propagation into unstable states.
\newblock {\em Physics Reports}, 386(2-6):29 -- 222, 2003.

\bibitem{MR2676976}
K.~Zumbrun.
\newblock A local greedy algorithm and higher-order extensions for global
  numerical continuation of analytically varying subspaces.
\newblock {\em Quart. Appl. Math.}, 68(3):557--561, 2010.

\bibitem{zumbrun98}
K.~Zumbrun and P.~Howard.
\newblock Pointwise semigroup methods and stability of viscous shock waves.
\newblock {\em Indiana Univ. Math. J.}, 47(3):741--871, 1998.

\end{thebibliography}

\end{document}